\documentclass[preprint]{elsarticle}
\usepackage[ansinew]{inputenc}
\usepackage{a4wide}
\usepackage{calc}
\usepackage{times}
\usepackage{amssymb,amsmath,amstext,amsthm,amsfonts, ams fonts}
\usepackage{colortbl,color,graphics,graphicx,epsfig}
\definecolor{gray}{rgb}{0.8,0.8,0.8}
\usepackage{dsfont}
\usepackage{bbm}
\usepackage{nicefrac}
\usepackage[pdflinkmargin=5pt,pdfstartview={FitBH -32768},plainpages=false]{hyperref}

\renewcommand{\P}{\mathbb{P}}
\newcommand{\F}{\mathcal{F}}
\renewcommand{\H}{\mathcal{H}}
\newcommand{\G}{\mathcal{G}}
\newcommand{\E}{\mathbb{E}}
\newcommand{\N}{\mathds{N}}
\newcommand{\R}{\mathds{R}}
\newcommand{\1}{\mathbbm{1}}
\renewcommand{\d}{\Delta}
\newcommand{\KLEINO}{{\scriptstyle{\mathcal{O}}}}
\newcommand{\oc}{\text{\textcircled{\texttt{1}}}}
\newcommand{\zc}{\text{\textcircled{\texttt{2}}}}
\newcommand{\dc}{\text{\textcircled{\texttt{3}}}}
\newcommand{\qc}{\text{\textcircled{\texttt{4}}}}

\DeclareMathAccent{\verywidehat}{\mathord}{largesymbols}{'144}
\newcommand{\var}{\mathbb{V}\hspace*{-0.05cm}\textnormal{a\hspace*{0.02cm}r}}

\newcommand{\cov}{\mathbb{C}\textnormal{o\hspace*{0.02cm}v}}
\DeclareMathOperator{\AVAR}{\mathbf{AVAR}}
\renewcommand{\:}{\mathrel{\mathop{:}}}

\newdefinition{remark}{Remark}
\newdefinition{defi}{Definition}
\newtheorem{theo}{Theorem}
\newtheorem{assump}{Assumption}
\newtheorem{prop}{Proposition}[section]
\newtheorem{lem}[prop]{Lemma}
\newtheorem{cor}[prop]{Corollary}

\begin{document}
\begin{frontmatter}



\title{An estimator for the quadratic covariation of asynchronously observed It\^{o} processes with noise: Asymptotic distribution theory}

\author{Markus Bibinger\fnref{label1}}
\fntext[label1]{Financial support from the Deutsche Forschungsgemeinschaft via SFB 649 `\"Okonomisches Risiko', Humboldt-Universität zu Berlin, is gratefully acknowledged.}

\address{Institut f\"ur Mathematik, Humboldt-Universit\"at zu Berlin, Unter den Linden 6, 10099 Berlin, Germany}

\begin{abstract}
The article is devoted to the nonparametric estimation of the quadratic covariation of non-synchronously observed It\^{o} processes in an additive microstructure noise model.
In a high-frequency setting, we aim at establishing an asymptotic distribution theory for a generalized multiscale estimator including a feasible central limit theorem with optimal convergence rate on convenient regularity assumptions. The inevitably remaining impact of asynchronous deterministic sampling schemes and noise corruption on the asymptotic distribution is precisely elucidated. A case study for various important examples, several generalizations of the model and an algorithm for the implementation warrant the utility of the estimation method in applications.

\end{abstract}
\begin{keyword}
non-synchronous observations \sep microstructure noise \sep integrated covolatility \sep multiscale estimator \sep stable limit theorem\\ \vspace*{.3cm}\noindent
\vspace*{.3cm}\noindent
\end{keyword}
\end{frontmatter}
\thispagestyle{plain}
\section{Introduction\label{sec:1}}
The nonparametric estimation of the univariate quadratic variation of a latent semimartingale from $n$ observations in a high-frequency setting with additive observation noise has been studied intensively in recent years. It is known from \cite{gloter} that $n^{\nicefrac{1}{4}}$ constitutes a lower bound for the rate of convergence. An important motivation which has stimulated an alliance of economists and statisticians to establish estimation techniques for this kind of latent semimartingale models is their utility in estimating daily integrated (co-)volatilities from high-frequency intraday returns that serve as a basis for risk management as well as portfolio optimization and hedging strategies.
The last years have seen an enormous increase of the amount of trading activities for many liquid securities. Paradoxically, the availability of high-frequency data necessitated a new angle on financial modeling. In fact, for every semimartingale the discrete realized (co-)volatilities converge in probability to the integrated measures. However, realized volatilities of typical high-frequency financial time series data explode for very high frequencies. This effect
is ascribed to market microstructure frictions. Sources of the market microstructure noise are manifold. One important role plays the occurrence of bid-ask spreads. Aside from that transaction costs, strategic trading, limited market depths and discreteness of prices spread out the structure of the long-run dynamics that can be characterized by semimartingales.\\
This strand of literature followed \cite{zhangmykland} that has attracted  a lot of attention to this estimation problem. The so-called two-scales realized volatility by \cite{zhangmykland} is based on subsampling and a bias-correction and a stable central limit theorem with $n^{\nicefrac{1}{6}}$-rate has been proved. A refinement of the subsample approach using multiple scales in \cite{zhang} and related alternative techniques in \cite{bn2}, \cite{preavg} and \cite{xiu} have led to rate-optimal estimators and feasible stable central limit theorems. For the more specific nonparametric model with Gaussian noise, asymptotic equivalence in the Le Cam sense to a Gaussian shift experiment is shown in \cite{reiss} and an asymptotically efficient estimator whose asymptotic variance equals the parametric efficiency bound is constructed.\\
In the article on hand we are concerned with a multivariate stetting and apart from taking additive microstructure noise into account, we we focus on a way to deal with non-synchronous observation schemes. This is also a central theme in financial applications. When realized covolatilities are calculated for fixed time distances and a previous-tick interpolation is applied, the phenomenon of the so-called Epps effect described in \cite{epps} appears that the realized covolatility tends to zero at the highest frequencies.\\
A methodology to deal with non-synchronous observations in a bivariate It\^{o} processes model has been proposed by \cite{hy}. The so-called Hayashi-Yoshida estimator has superseded simpler previous-tick interpolation methods setting the standard for the estimation of the quadratic covariation from asynchronous observations in the absence of microstructure noise effects.\\
Our estimation approach, first proposed in \cite{bibinger}, for the most general case in the presence of noise and non-synchronicity arises as a combination of the multiscale estimator to handle noise contamination on the one hand and a synchronization algorithm in accordance with the Hayashi-Yoshida estimator to cope with non-synchronicity on the other hand. A first attempt in the same direction, combining one-scale subsampling and the Hayashi-Yoshida estimator, has been given in \cite{palandri}.\\
In \cite{bibinger} it has been shown in the spirit of \cite{gloter} that the optimal convergence rate $n^{\nicefrac{1}{4}}$ carries over to the general multidimensional setup. The mathematical analysis of our generalized multiscale estimator in \cite{bibinger} shows that it is rate-optimal. \\
Alternative approaches to similar statistical models has been suggested by \cite{bn1}, \cite{kinnepoldi} and \cite{sahalia}. 
In \cite{bn1} a kernel-based method with a previous-tick interpolation to so-called refresh times is proposed and a stable central limit theorem with sub-optimal $n^{\nicefrac{1}{5}}$-rate is established for a multivariate non-synchronous design. This estimator, furthermore, ensures that the estimated covariance matrix is positive semi-definite. \cite{kinnepoldi} and \cite{sahalia} come up with combinations of pre-averaging (\cite{preavg},\cite{poldi}) and the Hayashi-Yoshida estimator and of the univariate quasi-maximum-likelihood method by \cite{xiu}, the polarization identity and a generalized synchronization scheme which is different from the Hayashi-Yoshida ansatz that we use, respectively, both also attaining the optimal rate.\\
In this article we aim at providing an asymptotic distribution theory for the generalized multiscale estimator. In distinction from alternative methods, the influence of non-synchronicity effects on the expectation is null and on the variance limited up to an interaction of interpolation steps and microstructure noise. The main result is a feasible stable central limit theorem for its estimation error with optimal rate and a closed-form asymptotic variance that does not hinge on interpolation errors in the signal term. The stable weak convergence of the estimation error to a centred mixed Gaussian limit and the consistent estimation of the random unknown asymptotic variances are the essential steps towards statistical inference and confidence sets. The theory is grounded on stable limit theorems for semimartingales from \cite{jacod1}.\\ 
In Section \ref{sec:2} we present the model and our main findings. Section \ref{sec:3} comes up with a concise overview on the construction of the estimator and in Section \ref{sec:4} we develop the asymptotic theory. In Section \ref{sec:5} we propose a consistent estimator for the asymptotic variance and Section \ref{sec:6} comprises various extensions and and a concluding discussion. The proofs are postponed to the Appendix.
\section{Model and key result\label{sec:2}}
The considered statistical model of noisy latently observed It\^{o} processes at deterministic observation times is precisely described by Assumptions \ref{eff}-\ref{e} in this section. 
\begin{assump}[\textbf{efficient processes}]\label{eff}
On a filtered probability space $\left(\Omega,\F,\left(\F_t\right),\P\right)$, the efficient processes $X=(X_t)_{t\in\R^{+}}$ and $Y=(Y_t)_{t\in\R^{+}}$ are It\^{o} processes defined by the following stochastic differential equations:
\begin{align*} dX_t&=\mu_t^X\,dt+\sigma_t^X\,dB_t^X~,\\
							dY_t&=\mu_t^Y\,dt+\sigma_t^Y\,dB_t^Y~,
\end{align*}
with two $\left(\F_t\right)$--adapted standard Brownian motions $B^X$ and $B^Y$ and $\rho_t\,dt=d\left[ B^X,B^Y\right]_t$.
The drift processes $\mu_t^X$ and $\mu_t^Y$ are $\left(\F_t\right)$--adapted locally bounded stochastic processes and the spot volatilities $\sigma_t^X$ and $\sigma_t^Y$ and $\rho_t$ are assumed to be $\left(\F_t\right)$--adapted with continuous paths. We assume strictly positive volatilities and the Novikov condition $\E\left[\exp{\left((1/2)\int_0^T(\mu^{\,\cdot\,}/\sigma^{\,\cdot\,})^2_t\,dt\right)}\right]<\infty$ for $X$ and $Y$.
\end{assump}
\begin{assump}[\textbf{observations}]\label{grid}
The deterministic observation schemes $\mathcal{T}^{X,n}=\{0\le t_0^{(n)}<t_1^{(n)}<\ldots<t_n^{(n)}\le T\}$ of $X$ and $\mathcal{T}^{Y,m}=\{0\le \tau_0^{(m)}<\tau_1^{(m)}<\ldots<\tau_m^{(m)}\le T\}$ of $Y$ are assumed to be regular in the following sense:
There exists a constant $0<\alpha\le 1/9$ such that 
\begin{subequations}
\begin{align}
\delta_n^X&=\sup_{i\in \{1,\ldots,n\}}{\left(\left(t_i^{(n)}-t_{i-1}^{(n)}\right),t_0^{(n)},T-t_n^{(n)}\right)}~\;\,=\mathcal{O}\left(n^{-\nicefrac{8}{9}-\alpha}\right)~,\\
\delta_m^Y&=\sup_{j\in \{1,\ldots,m\}}{\left(\left(\tau_j^{(m)}-\tau_{j-1}^{(m)}\right),\tau_0^{(m)},T-\tau_{m}^{(m)}\right)}=\mathcal{O}\left(m^{-\nicefrac{8}{9}-\alpha}\right)~.
\end{align}\end{subequations}
We consider asymptotics where the number of observations of $X$ and $Y$ are assumed to be of the same asymptotic order $n=\mathcal{O}(m)$ and $m=\mathcal{O}(n)$ and express that shortly by $n\sim m$.
The efficient processes $X$ and $Y$ which satisfy Assumption \ref{eff} are discretely observed at the times $\mathcal{T}^{X,n}$ and $\mathcal{T}^{Y,m}$ with additive observation noise:
\begin{subequations}
\begin{align}\tilde X_{t_i^{(n)}}=\int_0^{t_i^{(n)}}\mu_t^X\,dt+\int_0^{t_i^{(n)}}\sigma_t^X\,dB_t^X+\epsilon_{t_i^{(n)}}^X~,0\le i\le n~,\end{align}
\begin{align}\tilde Y_{\tau_j^{(m)}}=\int_0^{\tau_j^{(m)}}\mu_t^Y\,dt+\int_0^{\tau_j^{(m)}}\sigma_t^Y\,dB_t^Y+\epsilon_{\tau_j^{(m)}}^Y~,0\le j\le m~.\end{align}
\end{subequations}
\end{assump}
Although we consider sequences of deterministic observation times, the case of random sampling that is independent of the observed processes is included when regarding the conditional law.\\
It turns out that it is accurate to prove the key result of the article on the following i.\,i.\,d.\,assumption on the microstructure noise since a closed-form expression for the asymptotic variance is not available for a combination of general asynchronous observation schemes and serially dependent observation errors. Since an extension to non-i.\,i.\,d.\,noise is crucial for the utility in financial applications, we comment on the robustness of our estimator to that case in Section \ref{sec:6}. 
\begin{assump}[\textbf{microstructure noise}]\label{e}
The discrete microstructure noise processes $$\epsilon_{t_i^{(n)}}^X,\epsilon_{\tau_j^{(m)}}^Y,0\le i\le n,0\le j\le m~.$$ are centred i.\,i.\,d.\,, independent of each other and independent of the efficient processes $X$ and $Y$. We assume that the observation errors have finite fourth moments and denote the variances
\begin{align*}\eta_X^2=\var\left(\epsilon_{t_1^{(n)}}^X\right)~,~\eta_Y^2=\var\left(\epsilon_{\tau_1^{(m)}}^Y\right)~.\end{align*}
\end{assump}
The number of synchronized observations $N\sim n\sim m$ which appears in the rate of our feasible stable central limit theorem is introduced in Section \ref{sec:3}.
\begin{theo}[\textbf{feasible stable central limit theorem}]\label{fclt}
The generalized multiscale estimator \eqref{multiscale} specified by the later given weights \eqref{optimalweights}, with $M_N=c_{multi}\cdot \sqrt{N}$ converges on the Assumptions \ref{eff}, \ref{grid}, \ref{e} and further mild regularity conditions on the asymptotics of the sampling schemes, stated below in Assumptions \ref{grid3} and \ref{grid4}, $\mathcal{F}-$stably in law with optimal rate $N^{\nicefrac{1}{4}}\sim n^{\nicefrac{1}{4}}\sim m^{\nicefrac{1}{4}}$ to a mixed Gaussian limiting distribution:
\begin{align*}N^{\nicefrac{1}{4}}\left(\widehat{\left[ X,Y\right]}_T^{multi}-\left[ X,Y\right]_T\right)\stackrel{st}{\rightsquigarrow}\mathbf{N}\left(0,\AVAR_{multi}\right)\end{align*}
with a almost surely finite random asymptotic variance given in \eqref{avarmulti} in Theorem \ref{cltmulti}.
With the consistent estimator for the asymptotic variances $\widehat{\AVAR}_{multi}$ in Proposition \ref{avarest} , the feasible central limit theorem
\begin{align}n^{\nicefrac{1}{4}}\frac{\left(\widehat{\left[ X,Y\right]}_T^{multi}-\left[ X,Y\right]_T\right)}{\widehat{\AVAR}_{multi}}&\stackrel{st}{\rightsquigarrow}\mathbf{N}(0,1)~,\end{align}
holds true.
\end{theo}
The notion of stable weak convergence going back to \cite{renyi} is essential for our asymptotic theory. Stable weak convergence $X_n\stackrel{st}{\rightsquigarrow} X$ is the joint weak convergence of $(X_n,Z)$ to $(X,Z)$ for every measurable bounded random variable $Z$. The limiting random variables in stable limit theorems are defined on extensions of the original underlying probability spaces. The reason for us to involve this concept of a stronger mode of weak convergence is that mixed normal limiting distributions are derived where asymptotic variances are themselves strictly positive random variables. Provided we have  a consistent estimator $V_n^2$ for such a random asymptotic variance $V^2$ on hand, the stable central limit theorem $X_n\stackrel{st}{\rightsquigarrow}VZ$ with $Z$ distributed according to a standard Gaussian law, yields the joint weak convergence $(X_n,V_n^2)\rightsquigarrow (VZ,V^2)$ and also $X_n/V_n\rightsquigarrow Z$ and hence allows to perform statistical inference providing tests or confidence intervals.\\
In the proofs of our limit theorems we will `remove' the drifts in the sense that after a transformation to an equivalent martingale measure stable central limit theorems for It\^{o} processes without drift are proved and, as illustrated in \cite{inference}, stability of the weak convergence ensures that the asymptotic law holds true under the original measure. In this sense stable convergence is commutative with measure change.\\
From now on, we often omit the superscripts of observation times for a shorter notation.
\section{Brief review on the foundation\label{sec:3}}
\subsection{Subsampling and the multiscale estimator}
In the model imposed by Assumption \ref{eff}, Assumption \ref{grid} with synchronous observations, $n=m$ and $t_j^{(n)}=\tau_j^{(n)},1\le j\le n$, and Assumption \ref{e}, the realized (co-)volatilities do not provide consistent estimators for the quadratic (co-)variations any more. The variance due to noise conditional on the paths of the efficient processes
$$\var_{X,Y}\left(\sum_{j=1}^n\left(\tilde X_{t_j}-\tilde X_{t_{j-1}}\right)\left(\tilde Y_{t_j}-\tilde Y_{t_{j-1}}\right)\right)=4n\,\eta_X^2\eta_Y^2~,$$
increases linearly with $n$. The error due to noise perturbation can be reduced by the following estimator, which has been proposed for the univariate estimation of integrated volatility as the ``second best approach'' in \cite{zhangmykland} and which is called one-scale subsampling estimator in \cite{bibinger} and throughout this article.
\begin{figure}[t]\fbox{
\includegraphics[width=8.5cm]{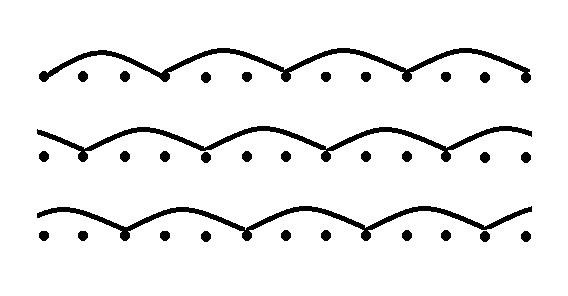}\includegraphics[width=6cm]{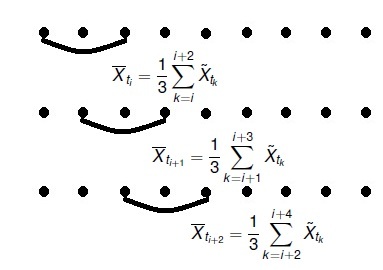}}
\caption{\label{subsample}Sketch of the subsampling approach.}
\end{figure}
It can be motivated from two perspectives that are both sketched in Figure \ref{subsample}. On the left-hand side we have visualized that one can calculate simultaneously lower frequent realized covolatilities using subsamples, e.\,g.\,to the lag three in Figure \ref{subsample}, and (post-)average them.
\begin{subequations}
\begin{align}\label{subsyn}\widehat{\left[ X,Y\right]}_T^{sub}=\frac{1}{i}\sum_{j=i}^n\left(\tilde X_{t_j}-\tilde X_{t_{j-i}}\right)\left(\tilde Y_{t_j}-\tilde Y_{t_{j-i}}\right)~.\end{align}
This motivation given in \cite{zhangmykland} is in line with the former common practice of a sparse-sampled low-frequency realized (co-)volatility estimator and proposes to use an average instead of one single lower frequent realized measure.\\
The same estimator arises as the usual realized covolatility calculated from the time series on that a linear filter is run before, what means that non-noisy observations at a time $t_j$ are estimated with a (pre-)average of noisy observations at times $t_j,\ldots,t_{j+i}$ for some $i$. This is sketched on the right-hand side of Figure \ref{subsample} for $i=3$. Passing over to increments leads to telescoping sums and we end up finally with the one-scale subsampling estimator.\\
Since on the Assumption \ref{e} there is no bias due to noise for the bivariate estimator, it already corresponds to the ``first best approach'' from \cite{zhangmykland} whereas in the univariate case a bias-correction completes the two scales realized volatility (TSRV):
\begin{align}\label{tsrv}\widehat{\left[ X\right]}_T^{TSRV}=\frac{1}{i}\sum_{j=i}^n\left(\tilde X_{t_j}-\tilde X_{t_{j-i}}\right)^2-\frac{1}{2n}\sum_{i=1}^n\left(\tilde X_{t_j}-\tilde X_{t_{j-1}}\right)^2~.\end{align}
\end{subequations}
There is a trade-off between the signal term and the error due to noise. Choosing $i=c_{sub}n^{\nicefrac{2}{3}}$ dependent on $n$ with a constant $c_{sub}$, the overall mean square error is minimized and of order $n^{-\nicefrac{1}{3}}$. The one-scale subsampling estimator \eqref{subsyn} is hence a consistent and asymptotically unbiased estimator. The rate of convergence $n^{\nicefrac{1}{6}}$, however, is slow and does not attain the optimal rate $n^{\nicefrac{1}{4}}$ determined in \cite{bibinger}.
For this reason, we focus on a multiscale extension of the subsampling approach on which the methods developed in \cite{bibinger} are based on. The multiscale realized covolatility (MSRC), and the univariate multiscale realized volatility (MSRV) introduced in \cite{zhang}, are linear combinations of one-scale subsampling estimators with $M_n$ different subsampling frequencies $i=1,\ldots,M_n$: 
\begin{subequations}
\begin{align}\label{multisyn}\widehat{\left[ X,Y\right]}_T^{multi}=\sum_{i=1}^{M_N}\frac{\alpha_{i,M_N}^{opt}}{i}\sum_{j=i}^n\left(\tilde X_{t_j}-\tilde X_{t_{j-i}}\right)\left(\tilde Y_{t_j}-\tilde Y_{t_{j-i}}\right)~,\end{align}
\begin{align}\label{msrv}\widehat{\left[ X\right]}_T^{multi}=\sum_{i=1}^{M_N}\frac{\alpha_{i,M_N}^{opt}}{i}\sum_{j=i}^n\left(\tilde X_{t_j}-\tilde X_{t_{j-i}}\right)^2~.\end{align}
\end{subequations}
The weights are chosen such that the estimator is asymptotically unbiased and the error due to noise minimized. They are given later in \eqref{optimalweights} and can be chosen equally for the bivariate and the univariate case. Those are the standard discrete weights of \cite{zhang} and we abstain from giving a more general class of possible weight functions.\\
The mean square error of the multiscale realized covolatility \eqref{multisyn} can be split in uncorrelated addends that stem from discretization, microstructure noise and cross terms and end-effects. They are of orders $M_n/n$, $n/M_n^3$, and $M_n^{-1}$, respectively. Hence, a choice $M_n=c_{multi}\sqrt{n}$ leads to a rate-optimal $n^{\nicefrac{1}{4}}$-consistent estimator.\\
The following stable central limit theorems for the multiscale realized covolatility \eqref{multisyn} and the one-scale estimator \eqref{subsyn} are implied by Theorem \ref{cltmulti} and Corollary \ref{cltone}:
\begin{prop}\label{cltsynchronous}
On Assumptions \ref{eff}, \ref{grid} and \ref{e} in the synchronous setup and if $(n/T)\sum_i(t_i^{(n)}-t_{i-1}^{(n)})^2$ converges to a continuously differentiable limiting function $G$ and the difference quotients converge uniformly to $G^{\prime}$ on $[0,T]$, the multiscale realized covolatility \eqref{multisyn} and the subsampling estimator \eqref{subsyn} converge stably in law to mixed normal limiting random variables:
\begin{subequations}
\begin{align}n^{\nicefrac{1}{4}}\left(\widehat{\left[ X,Y\right]}_T^{multi}-\left[ X,Y\right]_T\right)\stackrel{st}{\rightsquigarrow}\mathbf{N}\left(0\,,\,\AVAR_{multi,syn}\right)~,\end{align}
\begin{align}n^{\nicefrac{1}{6}}\left(\widehat{\left[ X,Y\right]}_T^{sub}-\left[ X,Y\right]_T\right)\stackrel{st}{\rightsquigarrow}\mathbf{N}\left(0\,,\,\AVAR_{sub,syn}\right)~,\end{align}
with
\begin{align}
\AVAR_{multi,syn}&=c_{multi}^{-3}24\eta_X^2\eta_Y^2+c_{multi}\frac{26}{35}T\int_0^TG^{\prime}(t)(\rho_t^2+1)(\sigma_t^X\sigma_t^Y)^2\,dt\\
&~\notag+c_{multi}^{-1}\frac{12}{5}\left(\eta_X^2\eta_Y^2+\eta_X^2\int_0^T(\sigma_t^Y)^2\,dt+\eta_Y^2\int_0^T(\sigma_t^X)^2\,dt\right)~,\end{align}
\begin{align}
\AVAR_{sub,syn}=c_{sub}^{-2}4\eta_X^2\eta_Y^2+c_{sub}\frac{2}{3}T\int_0^TG^{\prime}(t)(\rho_t^2+1)(\sigma_t^X\sigma_t^Y)^2\,dt~.
\end{align}
\end{subequations}
\end{prop}
\subsection{Synchronization and the Hayashi--Yoshida estimator\label{subsec:HY}}
We use the short notation $\d X_{t_i},i=1,\ldots,n$ from now on for increments $X_{t_i}-X_{t_{i-1}}$ and analogously for $Y$. The Hayashi-Yoshida estimator
\begin{align}\label{HYbase}\widehat{\left[  X, Y\right]}_T^{(HY)}=\sum_{i=1}^{n}\sum_{j=1}^{m}\d  X_{t_i}\d  Y_{\tau_j}\1_{[\min{(t_{i},\tau_{j})}>\max{(t_{i-1},\tau_{j-1})}]}~,\end{align}
where the product terms include all increments of the processes with overlapping observation time instants, has been proved in \cite{hy} to be consistent in a model of asynchronously observed It\^{o} processes with deterministic correlation, drift and volatility functions in the absence of observation noise and on further regularity conditions to be asymptotically normally distributed in \cite{hy2}.\\
For a combination of the strategy of the Hayashi-Yoshida estimator with techniques to handle noise conta\-mination, we use an iterative algorithm introduced in \cite{palandri} as `pseudo-aggregation'. Incorporating telescoping sums there are the following rewritings of the estimator \eqref{HYbase}:
\begin{align}\notag\widehat{\left[  X, Y\right]}_T^{(HY)}&=\sum_{i=1}^n\d X_{t_i}\left(Y_{t_{i,+}}-Y_{t_{i-1,-}}\right)\\
\label{HY} &=\sum_{i=1}^N\left( X_{g_i}- X_{l_i}\right)\left( Y_{\gamma_i}- Y_{\lambda_i}\right)\\
\notag  &=\sum_{i=1}^N\left( X_{T_{i,+}^X}- X_{T_{i-1,-}^X}\right)\left( Y_{T_{i,+}^Y}- Y_{T_{i-1,-}^Y}\right)~,
\end{align}
with the notion of next-tick interpolated times $t_{i,+}\:=\min_{0\le j\le m}{\left(\tau_j|\tau_j\ge t_i\right)}$ and previous-tick interpolated ones $t_{i,-}\:=\max_{0\le j\le m}{\left(\tau_j|\tau_j\le t_i\right)}$ in the first equality. This rewriting can be as well done in the symmetric way. \\
The illustration of \eqref{HY} that serves as a basis for the construction of the generalized multiscale estimator relies on an aggregation of the observations according to Algorithm \ref{A}. This algorithm, which is a concise version of the construction in \cite{bibinger}, stops after $(N+1)\le \min{(n,m)}+1$ steps when all observation times are grouped. Summation in \eqref{HY} can start with $i=0$ or $i=1$.
\addtocounter{figure}{-1}
\renewcommand{\figurename}{Algorithm}
\begin{figure}[ht]\fbox{
\begin{minipage}[c]{\textwidth} 
\renewcommand{\baselinestretch}{.75}\normalsize
Define ~~~~\,~$t^+(s)\:=\min{\{t_i\in \mathcal{T}^{X,n}|t_i\ge s\}}$, $\tau^+(s)\:=\min{\{\tau_j\in \mathcal{T}^{Y,n}|\tau_j\ge s\}}$;\\ \hspace*{1.5cm}$t^-(s)\:=\max{\{t_i\in \mathcal{T}^{X,n}|t_i< s\}}$, $\tau^-(s)\:=\max{\{\tau_j\in \mathcal{T}^{Y,n}|\tau_j< s\}}$.\\
first step:
\begin{itemize}\item For $t_{0}\le \tau_{0}$~~~$\dashrightarrow ~ g_0=t^+(\tau_0), l_0=t_0,\gamma_0=\tau_0, \lambda_0=\tau_0$.\\
\item For $t_{0}> \tau_{0}$~~~$\dashrightarrow ~ g_0=t_0, l_0=t_0,\gamma_0=\tau^+(t_0), \lambda_0=\tau_0$.
\end{itemize}
$i$th step (given $g_{i-1}$ and $\gamma_{i-1}$):
\begin{itemize}\item If $g_{i-1}=\gamma_{i-1}$ 
\begin{itemize} \item and $t^+(g_{i-1})\le \tau^+(\gamma_{i-1})$~~~$\dashrightarrow ~ g_i=t^+(\tau^+(\gamma_{i-1})), l_i=g_{i-1},\gamma_i=\tau^+(\gamma_{i-1}), \lambda_i=\gamma_{i-1}$.\\
\item and $t^+(g_{i-1})> \tau^+(\gamma_{i-1})$~~~$\dashrightarrow ~ g_i=t^+(g_{i-1}), l_i=g_{i-1},\gamma_i=\tau^+(t^+(g_{i-1})), \lambda_i=\gamma_{i-1}.$
\end{itemize}
\item If $g_{i-1}<\gamma_{i-1}$ 
\begin{itemize}
\item and $\gamma_{i-1}\le t^+(g_{i-1})$~~~$\dashrightarrow ~ g_i=t^+(g_{i-1}), l_i=g_{i-1},\gamma_i=\tau^+(t^+(g_{i-1})), \lambda_i=\tau^-(\gamma_{i-1})$.\\
\item and $\gamma_{i-1}>t^+(g_{i-1})$~~~$\dashrightarrow ~ g_i=t^+(\gamma_{i-1}), l_i=g_{i-1},\gamma_i=\gamma_{i-1}, \lambda_i=\tau^-(\gamma_{i-1})$.
\end{itemize}
\item If $g_{i-1}>\gamma_{i-1}$: symmetrically as for $g_{i-1}<\gamma_{i-1}$.
\end{itemize}
\end{minipage}}
\caption{\label{A}Iterative synchronization algorithm.}
\end{figure}
\renewcommand{\figurename}{Figure} 
\renewcommand{\baselinestretch}{1}\normalsize
\hspace*{-0.25cm}\noindent
In the last equality only the denotation expressions $g_i, \gamma_i, l_i, \lambda_i$ are substituted emphasizing that those sampling times obtained by Algorithm \ref{A} can be interpreted as previous- and next-tick interpolations again with respect to a synchronous sampling scheme $T_k\:=\min{(g_k,\gamma_k)},\,1\le k\le N$, which we call the closest synchronous approximation. Increments in \eqref{HY} are taken from previous-tick interpolations at left-end points of instants $[T_{k-1},T_k],\,2\le k\le N$ to next-tick interpolated sampling times at right-end points. Since $T_k=\max{(l_{k+1},\lambda_{k+1})},\,1\le k\le (N-1)$ holds true, we split the estimation error of \eqref{HY} in two uncorrelated parts $D_T^N+A_T^N$ with
\begin{align}D_T^N\:=\sum_{i=1}^N\left(\left(X_{T_i}-X_{T_{i-1}}\right)\left(Y_{T_i}-Y_{T_{i-1}}\right)-\int_{T_{i-1}}^{T_i}\rho_t\sigma_t^X\sigma_t^Y\,dt\right)\\
\notag~~~-\int_0^{t_0\wedge \tau_0}\rho_t\sigma_t^X\sigma_t^Y\,dt-\int_{t_n\wedge\tau_m}^T\rho_t\sigma_t^X\sigma_t^Y\,dt\end{align}
being the discretization error of a synchronous-type realized covolatility including in general non-observable idealized values at the times of the closest synchronous approximation and $A_T^N$ an additional error due to the lack of synchronicity, in particular next- and previous-tick interpolations. The times $T_k$ equal the so-called refresh times of \cite{bn1} and thus our synchronization differs from the one in \cite{bn1} by replacing pure previous-tick interpolation by the above given machinery of previous- and next-tick interpolations.\\
The asymptotic theory for the estimator \eqref{HY} as $N\rightarrow\infty$, concisely repeated here, is separately proved and presented in a more elaborate way in \cite{asyn}.\\
First, we take up the illustrative example from \cite{bibinger} to motivate the synchronization procedure. For further details and examples we refer to \cite{bibinger} and \cite{palandri}. 
Figure \ref{example2} visualizes the aggregation carried out by Algorithm \ref{A} and the times $T_i,\,i=0,\ldots,8$ for a toy example. The example emphasizes the important fact that consecutive right-end points of increments can be the same time points. The realized covolatility calculated from previous-tick interpolated values to refresh times equals
\begin{align*}(X_{t_2}-X_{t_0})(Y_{\tau_1}-Y_{\tau_0})+(X_{t_3}-X_{t_2})(Y_{\tau_3}-Y_{\tau_1})+(X_{t_5}-X_{t_3})(Y_{\tau_4}-Y_{\tau_3})+\\(X_{t_6}-X_{t_5})(Y_{\tau_5}-Y_{\tau_4})+(X_{t_7}-X_{t_6})(Y_{\tau_6}-Y_{\tau_5})+(X_{t_8}-X_{t_7})(Y_{\tau_7}-Y_{\tau_6})+\\(X_{t_9}-X_{t_8})(Y_{\tau_8}-Y_{\tau_7})+(X_{t_{10}}-X_{t_9})(Y_{\tau_{10}}-Y_{\tau_8})\end{align*}
and is systematically biased downwards by interpolations, whereas \eqref{HY} yields
\begin{align*}(X_{t_3}-X_{t_0})(Y_{\tau_1}-Y_{\tau_0})+(X_{t_3}-X_{t_2})(Y_{\tau_3}-Y_{\tau_1})+(X_{t_6}-X_{t_3})(Y_{\tau_4}-Y_{\tau_3})+\\(X_{t_7}-X_{t_5})(Y_{\tau_5}-Y_{\tau_4})+(X_{t_8}-X_{t_6})(Y_{\tau_6}-Y_{\tau_5})+(X_{t_8}-X_{t_7})(Y_{\tau_8}-Y_{\tau_6})+\\(X_{t_9}-X_{t_8})(Y_{\tau_9}-Y_{\tau_7})+(X_{t_{10}}-X_{t_9})(Y_{\tau_{10}}-Y_{\tau_8})~,\end{align*}
which is not biased due to interpolations.
\begin{figure}[t]
\begin{center}
\hspace*{-.2cm}\fbox{\includegraphics[width=14.4cm]{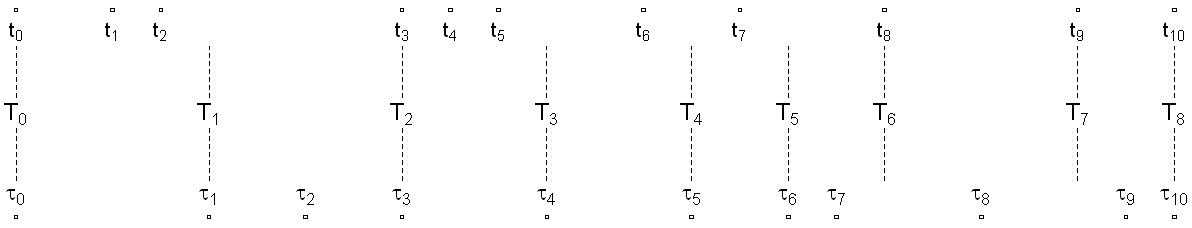}}
\end{center}\small
\begin{tabular}{||c|c|c|c|c|c|c|c|c|c||}
\hline
$i$&0&1&2&3&4&5&6&7&8\\
\hline
$\mathcal{H}^{i}$&$\{t_0\}$&$\{t_1,t_2,t_3\}$&$\{t_3\}$&$\{t_4,t_5,t_6\}$&$\{t_6,t_7\}$&$\{t_7,t_8\}$&$\{t_8\}$&$\{t_9\}$&$\{t_{10}\}$\\
$\mathcal{G}^{i}$&$\{\tau_0\}$&$\{\tau_1\}$&$\{\tau_2,\tau_3\}$&$\{\tau_4\}$&$\{\tau_5\}$&$\{\tau_6\}$&$\{\tau_7,\tau_8\}$&$\{\tau_8,\tau_9\}$&$\{\tau_9,\tau_{10}\}$\\
$l_i$&$t_0$&$t_0$&$t_2$&$t_3$&$t_5$&$t_6$&$t_7$&$t_8$&$t_9$\\
$g_i$&$t_0$&$t_3$&$t_3$&$t_6$&$t_7$&$t_8$&$t_8$&$t_9$&$t_{10}$\\
$\lambda_i$&$\tau_0$&$\tau_0$&$\tau_1$&$\tau_3$&$\tau_4$&$\tau_5$&$\tau_6$&$\tau_7$&$\tau_8$\\
$\gamma_i$&$\tau_0$&$\tau_1$&$\tau_3$&$\tau_4$&$\tau_5$&$\tau_6$&$\tau_8$&$\tau_9$&$\tau_{10}$\\
$T_i$&$\hspace*{-0.1cm}t_0\hspace*{-0.1cm}=\hspace*{-0.1cm}\tau_0\hspace*{-0.1cm}$&$\tau_1$&$\hspace*{-0.1cm}t_3\hspace*{-0.1cm}=\hspace*{-0.1cm}\tau_3$\hspace*{-0.1cm}&$\tau_4$&$\tau_5$&$\tau_6$&$t_8$&$t_9$&$\hspace*{-0.1cm}t_{10}\hspace*{-0.1cm}=\hspace*{-0.1cm}\tau_{10}\hspace*{-0.1cm}$\\
\hline
\end{tabular}\normalsize
\caption{\label{example2}Example for non-synchronous sampling design with the sets constructed by the synchronization algorithm, interpolated observation times occurring in the estimators and the synchronous approximation.}
\end{figure}
\begin{defi}[\textbf{quadratic (co-)variations of time}]\label{qvt}
For any $N\in\N$ let $T_i^{(N)},~i=0,\ldots,N$ be the times of the closest synchronous approximation and $g_i^{(N)},\gamma_i^{(N)},l_i^{(N)},\lambda_i^{(N)}$ the corresponding observation times that appear in the estimator \eqref{HY} defined above by Algorithm \ref{A} . $T/N$ is the mean of the time instants $\d T_i^{(N)}=T_i^{(N)}-T_{i-1}^{(N)},~i=1,\ldots,N$. Define the following sequences of functions
\begin{subequations}
\begin{align}\label{qvtg} G^N(t)=\frac{N}{T}\sum_{T_i^{(N)}\le t}\left(\d T_i^{(N)}\right)^2~,\end{align}
\begin{align}\label{qvtf}\notag F^N(t)=\frac{N}{T}\sum_{T_{i+1}^{(N)}\le t}(T_i^{(N)}-\lambda_i^{(N)})(g_i^{(N)}-T_i^{(N)})+\left(T_i^{(N)}-l_{i}^{(N)}\right)\left(\gamma_{i}^{(N)}-T_i^{(N)}\right)\\+\d T_{i+1}^{(N)}\left(T_i^{(N)}-l_{i+1}^{(N)}\right)+\d T_{i+1}^{(N)}\left(T_i^{(N)}-\lambda_{i+1}^{(N)}\right)~,\end{align}
\begin{align}\label{qvth} H^N(t)=\frac{N}{T}\sum_{T_{i+1}^{(N)}\le t}\left(T_i^{(N)}-l_{i+1}^{(N)}\right)\left(g_i^{(N)}-T_i^{(N)}\right)+\left(T_i^{(N)}-\lambda_{i+1}^{(N)}\right)\left(\gamma_i^{(N)}-T_i^{(N)}\right)~,\end{align}
\end{subequations}
for $t\in[0,T]$ that we call sequences of quadratic (co-)variations of times. 
\end{defi}
A stable central limit theorem for the estimation error is deduced in \cite{asyn} on the assumption that the sequences defined by \eqref{qvtg}, \eqref{qvtf} and \eqref{qvth} converge pointwise to continuous differentiable limiting functions $G,F,H$ and the sequences of difference quotients uniformly. The asymptotic quadratic variation of time $G$ of the $T_i^{(N)}$s influences the asymptotics of $D_T^N$. The covariation of times $F^N$ measures an interaction of interpolation errors between the two processes and $H^N$ the impact of the in general non-zero correlations of the products involving previous- and next-tick interpolations at the same $T_i^{(N)}$s for each process separately.\\
Consider as easiest example the synchronous equidistant sampling schemes with $N=n=m$ and $t_i^{(n)}=\tau_{j}^{(n)}=i/n,i=0,\ldots,n$. In this case $F^N$ and $H^N$ are identically zero since interpolations are redundant. The function $G^N$ is a step function that will tend to the identity on $[0,T]$ as $N\rightarrow\infty$.\\
Then, consider a situation of completely non-synchronous sampling schemes that originates from the complete synchronous equidistant one by shifting one time-scale half a time instant $1/2N$. We will call this situation intermeshed sampling. For this example, the synchronous approximation is still equidistant with instants $1/N$ and, hence, $G$ is the identity function. $F$ and $H$ are linear limiting functions with slope 1 and 1/4, respectively.\\
In \cite{asyn} we show for an important special case, independent homogeneous Poisson sampling, that the convergence assumptions on \eqref{qvtg}-\eqref{qvth} are fulfilled when replacing deterministic convergence by convergence in probability. Furthermore, the stochastic limits $G^{\prime}(t),F^{\prime}(t),H^{\prime}(t)$ are calculated explicitly. \\
The main result for the estimator \eqref{HY} is Theorem \ref{HYclt}. It serves as preparation to prove the stable limit theorem for the generalized multiscale estimator in Theorem \ref{cltmulti} and gives insight into the asymptotic distribution of \eqref{HY}. For the proof we refer to \cite{asyn}. A similar stable limit theorem for the original Hayashi-Yoshida estimator is provided in \cite{hy3}.
\begin{theo}\label{HYclt}The estimation error of \eqref{HY} converges on the Assumptions \ref{eff}, \ref{grid} and convergence assumptions on \eqref{qvtg}-\eqref{qvth} and the difference quotients stably in law to a centred, mixed Gaussian distribution:
\begin{equation}\sqrt{N}\left(\sum_{i=1}^N\left(X_{g_i}-X_{l_i}\right)\left(Y_{\gamma_i}-Y_{\lambda_i}\right)-\left[ X\,,\,Y\right]_T\right)\stackrel{st}{\rightsquigarrow}\mathbf{N}\left(0\,,\,v_{D_T}+v_{A_T}\right)~,\end{equation} 
with the asymptotic variance
\begin{equation*}\hspace*{-.05cm}v_{D_T}\hspace*{-.05cm}+\hspace*{-.05cm}v_{A_T}\hspace*{-.05cm}=\hspace*{-.05cm}T\hspace*{-.1cm}\int_0^T\hspace*{-.1cm} G^{\prime}(t)\hspace*{-.05cm}\left(\sigma_t^X\sigma_t^Y\right)^2\hspace*{-.05cm}\left(\rho_t^2+1\right)dt+T\hspace*{-.1cm}\int_0^T\hspace*{-.1cm}\left(F^{\prime}(t)\hspace*{-.05cm}\left(\sigma_t^X\sigma_t^Y\right)^2\hspace*{-.05cm}+2
H^{\prime}(t)\hspace*{-.05cm}\left(\rho_t\sigma_t^X\sigma_t^Y\right)^2\hspace*{-.05cm}\right)dt\end{equation*}
where the two addends come from the asymptotic variances of $D_T^N$ and $A_T^N$, respectively.
\end{theo}
\subsection{Hybrid approach to non-synchronous and noisy observations}
In \cite{bibinger} we have proposed the following combined estimation method for the quadratic covariation or integrated covolatility from noisy asynchronous observations. After applying Algorithm \ref{A} to the observation times, the generalized multiscale estimator is defined by
\begin{align}\label{multiscale}\widehat{\left[ X,Y\right]}_T^{multi}=\sum_{i=1}^{M_N}\frac{\alpha_{i,M_N}^{opt}}{i}\sum_{j=i}^N\left(\tilde X_{g_j^{(N)}}-\tilde X_{l_{j-i+1}^{(N)}}\right)\left(\tilde Y_{\gamma_j^{(N)}}-\tilde Y_{\lambda_{j-i+1}^{(N)}}\right)~.\end{align}
It is a weighted sum of $M_N$ one-scale subsampling estimators of the type
\vspace*{-.1cm}
\begin{align}\label{one-scale}\widehat{\left[ X,Y\right]}_T^{sub}=\frac{1}{i_N}\sum_{j=i_N}^N\left(\tilde X_{g_j^{(N)}}-\tilde X_{l_{j-i_N+1}^{(N)}}\right)\left(\tilde Y_{\gamma_j^{(N)}}-\tilde Y_{\lambda_{j-i_N+1}^{(N)}}\right)\end{align}
with subsampling frequencies $i=1,\ldots,M_N$ and optimal weights given later in \eqref{optimalweights}. Owing to the aggregation of non-synchronous observation times before applying subsampling and the multiscale approach, the resulting estimator has a conformable appearance as in the synchronous case \eqref{multisyn}. Recall that in the synchronous setting $g_j=\gamma_j=T_j$ and $l_{j-i+1}=\lambda_{j-i+1}=T_{j-i}$ holds.\\
Choosing $M_N=c_{multi} \cdot\sqrt{N}$ and $i_N=c_{sub}\cdot N^{\nicefrac{2}{3}}$, both estimators above provide consistent and asymptotically unbiased estimators with convergence rate $N^{\nicefrac{1}{4}}$ and $N^{\nicefrac{1}{6}}$, respectively. 
\section{Asymptotics and a stable central limit theorem for the generalized multiscale estimator\label{sec:4}}
A comprehensive analysis of the asymptotic distribution of the estimation error necessitates an elaborate screening of the conjunction of Algorithm \ref{A} and the joint sampling design $\left(\mathcal{T}^{X,n},\mathcal{T}^{Y,m}\right)$.\\
Note that the generalized multiscale estimator \eqref{multiscale} differs from the other plausible Hayashi-Yoshida version of a multiscale estimator
\begin{align}\label{alternativeestimator}\hspace*{-.15cm}\sum_{i=1}^{M_N}\frac{\beta_{i,M_N}^{opt}}{i}\hspace*{-.05cm}\sum_{j=i}^n\sum_{k\in\mathds{Z}}\hspace*{-.05cm}\Big(\tilde X_{t_j^{(n)}}\hspace*{-.05cm}-\tilde X_{t_{j-i}^{(n)}}\Big)\Big(\tilde Y_{\tau_{j+k\cdot i}^{(m)}}\hspace*{-.1cm}-\tilde Y_{\tau_{j+(k-1)\cdot i}^{(m)}}\Big)\hspace*{-.05cm}\1_{\{\max{(t_{j-i}^{(n)},\tau_{j+(k-1)\cdot i}^{(m)})}<\min{(t_j^{(n)},\tau_{j+k\cdot i}^{(m)})}\}}~,\end{align}
which arises as natural Hayashi-Yoshida multiscale estimator when, on the basis of (non-synchronized) observations of $\tilde X$ and $\tilde Y$, sparse-sample Hayashi-Yoshida estimators are averaged to one-scale subsample estimators and those extended to a linear combination using different time lags. We state without proof that this estimator is consistent, asymptotically unbiased and will attain the optimal rate of convergence. Nevertheless, we benefit from the data aggregation method and applying subsampling to the synchronized scheme, since the variance of our estimator \eqref{multiscale} is smaller than the one of this alternative estimator and we are able to find a feasible closed-form expression of the asymptotic variance.\\
The crucial difference between both approaches is that for the alternative method next- and previous-tick interpolation errors take place on sparse-sampling time intervals in average of order $i/N$ whereas the interpolation errors of the generalized multiscale estimator \eqref{multiscale} take place on the highest-frequency-scale and hence on intervals in average of order $1/N$.
In particular the decomposition
$$\left(X_{g_j^{(N)}}- X_{l_{j-i_N+1}^{(N)}}\right)=\big(\underbrace{ X_{g_j^{(N)}}- X_{{T_j}^{(N)}}}_{=\mathcal{O}_p(N^{-1/2})}+\underbrace{ X_{{T_j}^{(N)}}- X_{T_{j-i}^{(N)}}}_{=\mathcal{O}_p((i/N)^{(1/2)})}+\underbrace{ X_{T_{j-i}^{(N)}}-X_{l_{j-i_N+1}^{(N)}}}_{=\mathcal{O}_p(N^{-1/2})}\big)$$
of the increments of $X$ and analogously for $Y$, give an heuristic that the interpolation errors driving the error due to non-synchronicity asymptotically not affect the variance of the signal term. The stochastic orders are given for times instants of average order $N^{-1}$.\\
For a rigorous clarification of the asymptotic error due to noise and the cross term, both influenced by the i.\,i.\,d.\,observation errors at times $g_i,l_i,\gamma_i,\lambda_i$, we figure out the times $g_i=g_{i+1}$ and the right-end points $g_i^{(N)}=l_{i+1}^{(N)}, g_i^{(N)}=l_{i+2}^{(N)}$ that are as well preceding left-end points and analogously for the sampling times of $\tilde Y$. \\
All observation times $\gamma_i,\lambda_i$ are characterized through one of the following four mutually exclusive cases. Denote $\gamma_{j,-}$ the last observation time of $\tilde Y$ before $\gamma_j$ and $\gamma_{j,+}$ the first one after $\gamma_j$. We illustrate the allocation of the observation times for $\mathcal{T}^{Y,m}$ and $\gamma_j\,,\,j=1,\ldots,N-2$:
\vspace*{.1cm}
\begin{align*}&\oc ~\gamma_j\le g_j \hspace*{4.835cm}\Rightarrow \gamma_j \ne \gamma_{j+1}~,~\gamma_j=\lambda_{j+1}~,~\gamma_j\ne \lambda_{j+2}~~,\\
							&\zc~ \gamma_j> g_j~, ~\gamma_j\ge g_{j,+} \hspace*{3.1cm}\Rightarrow \gamma_j = \gamma_{j+1}~,~\gamma_j\ne\lambda_{j+1}~,~\gamma_j= \lambda_{j+2}~~,\\
							&\dc ~\gamma_j> g_j~, ~\gamma_j< g_{j,+}~,~\gamma_{j,+}>g_{j,+} \hspace*{1.1cm}\Rightarrow \gamma_j \ne \gamma_{j+1}~,~\gamma_j\ne\lambda_{j+1}~,~\gamma_j=\lambda_{j+2}~~, \vspace*{.1cm}\\
								&\qc~ \gamma_j> g_j~, ~\gamma_j< g_{j,+}~,~\gamma_{j,+}\le g_{j,+} \hspace*{1.1cm}\Rightarrow \gamma_j \ne \gamma_{j+1}~,~\gamma_j\ne\lambda_{j+1}~,~\gamma_j\ne\lambda_{j+2}~,~\gamma_{j,+}=\lambda_{j+2}~~.
								\end{align*}
Only sampling times distributed to case $\zc$ lead to repeated $\gamma_i=\gamma_{i+1}$. In cases $\oc,\zc$ and $\dc$ a subsequent left-end point $\lambda_k,k=i+1$ or $k=i+2$ of observation time instants incorporated in the subsampling estimators is designated by $\gamma_i$. All other $\lambda_k,\,k=2,\ldots,N$ appear in an allocation of sampling times of the type $\qc$, where $\lambda_{j+2}=\gamma_{j,+}\ne\gamma_l \,\forall l$. Recall that $\lambda_i\ne \lambda_k$ for all $i\ne k$ holds true.\\
If $\zc$ holds for $\gamma_j$ with fixed $j\in\{1,\ldots,N-2\}$ and if $k\:=\arg\min_{k\in\{j,\ldots,N-1\}}$${\left(\gamma_k>g_k\,,\,\gamma_k\ge g_{k,+}\right)}$ exists, then $\zc$ holds necessarily for one $g_l,l\in\{j+1,\ldots,k-1\}$ or $g_l=\gamma_l$.\\
In Table \ref{tabcases} we list the relations for the sampling design of our previous example.
\begin{table}[t]
\begin{center}
\begin{tabular}{|c|cccccccc|}
\hline
$i$&1&2&3&4&5&6&7&8\\
case $X$&$\zc$&$\oc$&$\dc$&$\dc$&$\zc$&$\oc$&$\oc$&$\oc$\\
case $Y$&$\oc$&$\oc$&$\oc$&$\oc$&$\oc$&$\dc$&$\qc$&$\oc$\\
relations $X$&$g_1=g_2=l_3$&$g_2=l_3$&$g_3=l_5$&$g_4=l_6$&$g_5=l_7$&$g_6=l_7$&$g_7=l_8$&--\\
relations $Y$&$\gamma_1=\lambda_2$&$\gamma_2=\lambda_3$&$\gamma_3=\lambda_4$&$\gamma_4=\lambda_5$&$\gamma_5=\lambda_6$&$\gamma_6=\lambda_8$&--&--\\
\hline
\end{tabular}
\end{center}
\caption{\label{tabcases}Allocation of sampling times to cases $\oc-\qc$ for the example.}
\end{table}
\begin{assump}[\textbf{asymptotic quadratic variation of time}]\label{grid3}
Assume that for the sequences of sampling schemes and the times $T_i^{(N)}$ of the closest synchronous approximations and for the sequence of quadratic variations of time $G^N(t)$ defined in Definition \ref{qvt}, the following holds true:
\begin{enumerate}
\item[(i)]
$G^N(t)\rightarrow G(t)$ as $N\rightarrow \infty$, where $G(t)$ is a continuously differentiable function on $[0,T]$.
\item[(ii)]For any null sequence $(h_N),\,h_N=\mathcal{O}\left(N^{-1}\right)$
\begin{align}\label{aqvtg2}\frac{G^N(t+h_N)-G^N(t)}{h_N}\rightarrow G^{\prime}(t)\end{align} 
uniformly on $[0,T]$ as $N\rightarrow \infty$..
\item[(iii)]The derivative $G^{\prime}(t)$ is bounded away from zero.
\end{enumerate}
\end{assump}
\begin{defi}[\textbf{degree of regularity of asynchronicity}]\label{dr}
For $N\in \N$ and sets $\mathcal{H}^{i},\mathcal{G}^{i},i=0,\ldots,N$ constructed from aggregated sampling schemes $\mathcal{T}^{X,n},\mathcal{T}^{Y,m}$ that fulfill Assumption \ref{grid}, define the following sequences of functions:
\begin{subequations}
\begin{align}I_X^N(t)=\frac{1}{N}\sum_{g_j^{(N)}\le t}\1_{\{g_j^{(N)}=g_{j-1}^{(N)}\}}~,\end{align}
\begin{align}I_Y^N(t)=\frac{1}{N}\sum_{\gamma_j^{(N)}\le t}\1_{\{\gamma_j^{(N)}=\gamma_{j-1}^{(N)}\}}~,\end{align}
\end{subequations}
which describe the degree of regularity of asynchronicity between observation times $\mathcal{T}^{X,n}$ and $\mathcal{T}^{Y,m}$.
\end{defi}
In the completely asynchronous case, we can directly conclude that $|I^N_X(t)-I^N_Y(t)|\le T/N$ for all $t\in[0,T]$ and one sequence suffices to reflect the regularity of the non-synchronous sampling schemes.
\begin{assump}[\textbf{asymptotic degree of regularity of asynchronicity}]\label{grid4}
Assume that for the sequences of sampling schemes and for the sequences of functions $I_X^N,I^N_Y$ defined in Definition \ref{dr}, the following holds true:
\begin{enumerate}
\item[(i)]
$I^N_X(t)\rightarrow I_X(t),I^N_Y(t)\rightarrow I_Y(t)$ as $N\rightarrow \infty$, where $I_X(t),I_Y(t)$ are continuously differentiable functions on $[0,T]$.
\item[(ii)]For any null sequence $(h_N),\,h_N=\mathcal{O}\left(N^{-1}\right)$
\begin{subequations} \begin{align}\frac{I_X^N(t+h_N)-I_X^N(t)}{h_N}\rightarrow I_X^{\prime}(t)~,\end{align} 
\begin{align}\label{adr}\frac{I_Y^N(t+h_N)-I_Y^N(t)}{h_N}\rightarrow I_Y^{\prime}(t)\end{align}\end{subequations}
uniformly on $[0,T]$ as $N\rightarrow \infty$.
\end{enumerate}
\end{assump}
For both, synchronous and intermeshed sampling which have been introduced in the last section, the sequences of functions $I_X^N,I_Y^N$ are identically zero. The functions defined in Definition \ref{dr} are non-negative and bounded above by $1$. In Section \ref{sec:6} we explicitly deduce the asymptotic degree of regularity of asynchronicity for mutually independent homogeneous Poisson sampling schemes. The term (asymptotic) degree of regularity of asynchronicity has been chosen since Assumption \ref{grid4} holds for all non-degenerate sequences where observation times conforming to one of the cases $\oc-\qc$ from above tend to be distributed according to some regular pattern and it gives information on the interaction of allocations of observation times. \\
It is interesting and might seem surprising at first glance that the asymptotics of the estimator \eqref{multiscale} hinges on this asymptotic feature whereas, as indicated before, the asymptotic interpolations to the closest synchronous approximation are asymptotically immaterial. This circumstance is caused by the fact that for the construction of an estimator with Algorithm \ref{A}, as for the original Hayashi-Yoshida estimator \eqref{HY}, observed values of the processes at next-tick interpolated observation times can appear twice. If there is observation noise, the number of observations allocated conforming to case $\zc$ has an impact on the asymptotics. The influence of interpolations is asymptotically vanishing for the combined method in contrast to the estimator \eqref{HY} with faster convergence rate $\sqrt{N}$ since interpolation steps take place on the time-scale of high-frequency observations, but lower-frequency sparse-sampled increments of the synchronous approximation are involved to reduce the error due to noise. 
We continue with the central result of this article:
\begin{theo}[\textbf{Central limit theorem for the generalized multiscale estimator}]\label{cltmulti}
On the Assumptions \ref{eff}, \ref{grid}, \ref{e}, \ref{grid3} and \ref{grid4}, the generalized multiscale estimator \eqref{multiscale} with noise-optimal weights $\alpha_{i,M_N}^{opt}=(12i^2/M_N^3)-(6i/M_N^2)\left(1+\KLEINO(1)\right)$, that are explicitly given in \eqref{optimalweights}, and $M_N=c_{multi}\cdot \sqrt{N}$ converges $\mathcal{F}-$stably in law with optimal rate $N^{\nicefrac{1}{4}}$ to a mixed Gaussian limiting distribution:
\begin{align*}N^{\nicefrac{1}{4}}\left(\widehat{\left[ X,Y\right]}_T^{multi}-\left[ X,Y\right]_T\right)\stackrel{st}{\rightsquigarrow}\mathbf{N}\left(0,\AVAR_{multi}\right)\end{align*}
with the asymptotic variance
\begin{align}\notag \AVAR_{multi}&=c_{multi}^{-3}\,\underbrace{\left(24+12\,(I_X(T)+I_Y(T))\right)\eta_X^2\eta_Y^2}_{=\AVAR_{\text{noise}}}+c_{multi}^{-1}\,\frac{12\eta_X^2\eta_Y^2}{5}\\ \label{avarmulti}
&\; +c_{multi}\,\underbrace{\frac{26}{35}T\int_0^TG^{\prime}(t)(\sigma_t^X\sigma_t^Y)^2(1+\rho_t^2)\,dt}_{=\AVAR_{dis,multi}}\\
\notag &\; +c_{multi}^{-1}\,\underbrace{\frac{12}{5}\left(\eta_Y^2\int_0^T(1+I^{\prime}_Y(t))(\sigma_t^X)^2\,dt\,+\eta_X^2\int_0^T(1+I^{\prime}_X(t))(\sigma_t^Y)^2\,dt\right)}_{=\AVAR_{cross}}~.\end{align}
\end{theo}
The weak convergence is proved to be stable with respect to the $\sigma$-algebra $\mathcal{F}$ associated with the efficient processes.
As a side result, we also obtain a stable central limit theorem for a simpler one-scale subsampling estimator:
\begin{cor}[\textbf{Central limit theorem for the one-scale subsampling estimator}]\label{cltone}
On the Assumptions \ref{eff}, \ref{grid}, \ref{e} and \ref{grid3}, the one-scale subsampling estimator with subsampling frequency $i_N=c_{sub}\cdot N^{\nicefrac{2}{3}}$ converges $\mathcal{F}$-stably in law with rate $N^{\nicefrac{1}{6}}$ to a mixed Gaussian limiting distribution:
\begin{align}N^{\nicefrac{1}{6}}\left(\widehat{\left[ X,Y\right]}_T^{sub}-\left[ X,Y\right]_T\right)\stackrel{st}{\rightsquigarrow}\mathbf{N}\left(0,\AVAR_{sub}\right)~,\end{align}
with the asymptotic variance
\begin{align}\label{avarone}\AVAR_{sub}&=c_{sub}^{-2}\,\underbrace{4\eta_X^2\eta_Y^2}_{=\AVAR_{\text{noise,sub}}}\hspace*{-.35cm}+\,c_{sub}\,\underbrace{\frac{2}{3}T\int_0^TG^{\prime}(t)(\sigma_t^X\sigma_t^Y)^2(1+\rho_t^2)\,dt}_{=\AVAR_{dis,sub}}
~.\end{align}
\end{cor}
For the proof of Theorem \ref{cltmulti}, we split the total estimation error of the generalized multiscale estimator in three asymptotically uncorrelated addends due to noise, cross terms and the signal term. For the one-scale subsampling estimator we follow the same ansatz. 
The orders of the errors have been derived in \cite{bibinger} and we focus on the asymptotic distribution here. \\
The error due to microstructure noise of the one-scale subsampling estimator has expectation zero and the variance yields
$$i_N^{-2}\sum_{j=i_N}^N\E\left[\left(\epsilon_{g_j}^X-\epsilon_{l_{j-i_N+1}}^X\right)^2\left(\epsilon_{\gamma_j}^Y-\epsilon_{\lambda_{j-i_N+1}}^Y\right)^2\right]=4Ni_N^{-2}\eta_X^2\eta_Y^2+\KLEINO\left(Ni_N^{-2}\right)~,$$
since observation noises of $\tilde X$ and $\tilde Y$ are independent of each other by Assumption \ref{e} and $l_k\ne l_r$ for $k\ne r$, $\lambda_k\ne \lambda_r$ for $k\ne r$ and if $g_k=g_{k+1}\Rightarrow \gamma_k<\gamma_{k+1}~,0\le k\le (N_1)$. Hence, the error due to noise is a sum of uncorrelated centred random variables with equal variances and the standard central limit theorem applies.\\
For the generalized multiscale estimator, we further decompose the error due to noise in a main part of order $N^{\nicefrac{1}{2}}M_N^{-\nicefrac{3}{2}}$ and two terms due to end-effects of orders $M_N^{-\nicefrac{1}{2}}$, where all three terms are asymptotically uncorrelated. In Propositions \ref{propnoise1} and \ref{propnoise2} we prove central limit theorems for these terms. Asymptotic normality holds both, conditionally and unconditionally on the paths of the efficient processes.\\
The error due to noise of the one-scale estimator does not depend on any further influence of the sampling schemes except the number of constructed sets $N$ and $G^{\prime}$. Cross terms are, in contrast to the multiscale case, asymptotically negligible since
\begin{align*}
&\E\left[\left(i_N^{-1}\sum_{j=i_N}^N\left(\left(X_{g_j}-X_{l_{j-i_N+1}}\right)\left(\epsilon_{\gamma_j}^Y-\epsilon_{\lambda_{j-i_N+1}}^Y\right)+\left(Y_{\gamma_j}-Y_{\lambda_{j-i_N+1}}\right)\left(\epsilon_{g_j}^X-\epsilon_{l_{j-i_N+1}}^X\right)\right)\right)^2\right]\\
&=i_N^{-2}\sum_{j=i_N}^N\left(\int_{l_{j-i_N+1}}^{g_{j}}(\sigma_t^X)^2\,dt~2\eta_Y^2+\int_{\lambda_{j-i_N+1}}^{\gamma_j}(\sigma_t^Y)^2\,dt~2\eta_X^2\right)+\KLEINO\left(i_N^{-2}\right)\\
&=\mathcal{O}\left(i_N^{-2}\right)=\KLEINO(1)~.\end{align*}
For the generalized multiscale estimator instead the cross terms are of order $M_N^{-\nicefrac{1}{2}}$ and will have effect upon the asymptotic distribution. In Proposition \ref{cltcross} a limit theorem is stated where the weak convergence also holds conditionally and unconditionally on the paths of the efficient processes. The asymptotic variance $\AVAR_{cross}$ includes the influence of the asymptotic degree of regularity of asynchronicity.\\
The error due to discretization of the one-scale subsampling estimator yields:
\vspace*{-.25cm}{\allowdisplaybreaks[3]{
\begin{align*}
&\frac{1}{i}\sum_{j=i}^N\left(X_{g_j}-X_{l_{j-i+1}}\right)\left(Y_{\gamma_j}-Y_{\lambda_{j-i+1}}\right)-\left[ X,Y\right]_T\\
&=\frac{1}{i}\sum_{j=i}^N\left(X_{T_j}-X_{T_{j-i}}\right)\left(Y_{T_j}-Y_{T_{j-i}}\right)-\left[ X,Y\right]_T\\
&+\frac{1}{i}\sum_{j=i}^N\left[\left(X_{g_j}-X_{T_j}\right)\left(Y_{T_j}-Y_{T_{j-i}}\right)+\left(X_{g_j}-X_{T_j}\right)\left(Y_{T_{j-i}}-Y_{\lambda_{j-i+1}}\right)\right. \\
&\left.~~~~~~ +\left(X_{T_{j-i}}-X_{l_{j-i+1}}\right)\left(Y_{\gamma_j}-Y_{T_{j}}\right)+\left(X_{T_{j-i}}-X_{l_{j-i+1}}\right)\left(Y_{T_j}-Y_{T_{j-i}}\right)\right.\\
&\left.~~~~~~ +\left(Y_{\gamma_j}-Y_{T_{j}}\right)\left(X_{T_j}-X_{T_{j-i}}\right)+\left(Y_{T_{j-i}}-Y_{\lambda_{j-i+1}}\right)\left(X_{T_j}-X_{T_{j-i}}\right)\right] \displaybreak[0]\\
&= \frac{1}{i}\sum_{j=i}^N\Big(\sum_{k=j-i+1}^j\Delta X_{T_k}\Big)\Big(\sum_{k=j-i+1}^j\Delta Y_{T_k}\Big)-\left[ X,Y\right]_T\\
&  +\frac{1}{i}\sum_{j=i}^{N-1}\left[\left(X_{g_j}-X_{T_j}\right)\Big(\sum_{k=j-i+1}^j\Delta Y_{T_k}\Big)+\left(X_{g_j}-X_{T_j}\right)\left(Y_{T_{j-i}}-Y_{\lambda_{j-i+1}}\right)\right.\\
&\left.~~~~~~ +\left(X_{T_{j-i}}-X_{l_{j-i+1}}\right)\left(Y_{\gamma_j}-Y_{T_{j}}\right)+\left(X_{T_{j}}-X_{l_{j+1}}\right)\Big(\sum_{k=j-i+1}^j\Delta Y_{T_k}\Big)\right.\\
&\left.~~~~~~ +\left(Y_{\gamma_j}-Y_{T_{j}}\right)\Big(\sum_{k=j-i+1}^j\Delta X_{T_k}\Big)+\left(Y_{T_{j}}-Y_{\lambda_{j+1}}\right)\Big(\sum_{k=j-i+1}^j\Delta X_{T_k}\Big)\right]+\mathcal{O}_p\left(N^{-1}\right) \displaybreak[0]\\
&=\frac{1}{i}\sum_{j=i}^N\Big(\sum_{k=j-i+1}^j\Delta X_{T_k}\Delta Y_{T_k}+\sum_{\substack{l\ne r \\ l,r\in\{j-i+1,\ldots,j\}}}\Delta X_{T_l}\Delta Y_{T_r}\Big)-\left[ X,Y\right]_T\\
&+\frac{1}{i}\sum_{j=i}^{N-1}\left[\left(X_{g_j}-X_{T_j}\right)\Big(\sum_{k=j-i+1}^j\Delta Y_{T_k}\Big)+\left(Y_{\gamma_j}-Y_{T_{j}}\right)\Big(\sum_{k=j-i+1}^j\Delta X_{T_k}\Big)\right.\\
& \left.~~~~~~+ \Delta X_{T_{j+1}}\Big(\sum_{k=j-i+1}^j\left(Y_{T_k}-Y_{\lambda_{k+1}}\right)\Big)+\Delta Y_{T_{j+1}}\Big(\sum_{k=j-i+1}^j\left(X_{T_k}-X_{l_{k+1}}\right)\Big)\right]\\
&~~~~~~+\mathcal{O}_p\left(i^{-1}N^{-\frac{1}{2}}\right)+\mathcal{O}_p\left(N^{-1}\right)~.\end{align*}}}
We have written the overall discretization error of the one-scale estimator as the sum of a discretization error of the closest synchronous approximation
\begin{align}\label{diserr}\hspace*{-.05cm}\sum_{j=1}^N\left(\hspace*{-.05cm}\Delta X_{T_j}\sum_{l=1}^{i\wedge j}\left(1-\frac{l}{i}\right)\Delta Y_{T_{j-l}}\hspace*{-.05cm}+\Delta Y_{T_j}\sum_{l=1}^{i\wedge j}\left(1-\frac{l}{i}\hspace*{-.05cm}\right)\Delta X_{T_{j-l}}\right)\hspace*{-.05cm}+\mathcal{O}_p\left(iN^{-1}\right)\hspace*{-.05cm}+\mathcal{O}_p\left(N^{-\nicefrac{1}{2}}\right)\end{align}
and the asymptotically negligible error due to the lack of synchronicity.
A stable central limit theorem using the theory of \cite{jacod1} for the leading term of order $i^{\nicefrac{1}{2}}N^{-\nicefrac{1}{2}}$ that will drive the asymptotic distribution is postponed to Proposition \ref{cltmsdis}. The error due to asynchronicity is treated in Proposition \ref{Asy}.\\
The discretization error of the generalized multiscale estimator is of order $M_N^{\nicefrac{1}{2}}N^{-\nicefrac{1}{2}}$ and that of the one-scale estimator of order $i_N^{\nicefrac{1}{2}}N^{-\nicefrac{1}{2}}$. There is a trade-off between the error due to noise and the discretization error for both estimators. For the generalized multiscale estimator these are of orders $N^{\nicefrac{1}{2}}M_N^{-\nicefrac{3}{2}}$ and $M_N^{\nicefrac{1}{2}}N^{-\nicefrac{1}{2}}$, respectively. Remaining other terms are of orders $M_N^{-\nicefrac{1}{2}}$. Thus, choosing $M_N=c_{multi}\cdot N^{\nicefrac{1}{2}}$, the total estimation error is minimized and of order $M_N^{-\nicefrac{1}{2}}=N^{-\nicefrac{1}{4}}$ which constitutes the optimal rate of convergence in Theorem \ref{cltmulti}.\\ The weak convergence of the discretization error is proved to be stable, so it converges jointly in law with every bounded $\mathcal{F}$-measurable random variable defined on the same probability space. Since the asymptotic normality of the cross term and the error due to noise holds both, conditionally and unconditionally given the efficient processes, and the discretization error is independent of $\epsilon^X$ and $\epsilon^Y$ we can apply a central limit theorem for mixing triangular arrays as in \cite{utev} to the sum that is adapted with respect to $\mathcal{A}_j=\sigma\left(\epsilon_{t_k}^X,t_k<T_{j+1},\epsilon_{\tau_k}^Y,\tau_k<T_{j+1},\mathcal{F}_{T_j}\right)$ where $\mathcal{F}$ is the $\sigma$-algebra associated with the efficient processes. The asymptotic variance is the sum of those of the uncorrelated addends. With the Cram\'{e}r-Wold device joint normality and asymptotic independence of the different errors can be concluded.\\
This is likewise for the one-scale estimator and Corollary \ref{cltone}. Choosing the subsampling frequency $i_N=c_{sub}\cdot N^{\nicefrac{2}{3}}$ balances the variance of the error due to noise which is of order $N i^{-2}$ and the discretization variance of order $iN^{-1}$. 
\section{Asymptotic variance estimation\label{sec:5}}
The asymptotic variances \eqref{avarmulti} and \eqref{avarone} of the generalized multiscale estimator \eqref{multiscale} and the one-scale subsampling estimator
\eqref{one-scale}, appearing in the stable central limit theorems in Theorem \ref{cltmulti} and Corollary \ref{cltone}, are random and depend on unknown quantities. In this section, we aim at estimating these asymptotic variances consistently to make our limit theorems feasible.\\
It is a known result that a consistent estimator of the noise variances is given by \eqref{noisevarest} (cf.\,\cite{zhangmykland}). Furthermore, the estimators for $\eta_X^2$ and $\eta_Y^2$ are asymptotically uncorrelated on Assumption \ref{e}, since the uncorrelated noise terms dominate the correlated Brownian parts. The constant $I_X(T)+I_Y(T)$ in the noise part of \eqref{avarmulti} can be estimated with the empirical version $I_X^N(T)+I_Y^N(T)$ that converges as $N\rightarrow\infty$ on Assumption \ref{grid4}. Eventually, consistent estimators for the discretization variances and the variance due to cross terms for the multiscale estimator are required.\\
We propose histogram-type estimators using bins according to timescales associated with the quadratic variation of synchronized sampling times and associated with the degree of regularity of asynchronicity, respectively. For this purpose, given a chosen number of bins $K_N$, with $K_N\rightarrow\infty$ and $K_N^{-1}N\rightarrow\infty$ as $N\rightarrow\infty$, we determine the assigned non-equispaced  bin-widths $\Delta G_j^N=G_j^N-G_{j-1}^N$, $\Delta {(I_X)}_j^N={(I_X)}_j^N-{(I_X)}_{j-1}^N$ and $\Delta {(I_Y)}_j^N={(I_Y)}_j^N-{(I_Y)}_{j-1}^N$, $j\in\{1,\ldots,K_N\}$, where
\begin{align*}
G_j^N\:=\inf{\big\{t\in[0,T]\,\big|\;G^N(t)=(N/T)\sum_{T_k^{(N)}\le t}\big(\Delta T_k^{(N)}\big)^2\ge (j/{K_N})\cdot G^N(T)\big\}}~,\end{align*}
$j\in\{1,\ldots,K_N\}$, and analogously for the functions $I_X^N$ and $I_Y^N$ if $I_X^N(T)>0$ and $I_Y^N(T)>0$. Set $G_0^N=(I_X)_0^N=(I_Y)_0^N\:=0$ and recall that those functions are monotone increasing on $[0,T]$. On each bin we calculate multiscale estimators in the same spirit as \eqref{multiscale} and its univariate version from \cite{zhang} for the increase of the quadratic (co-) variations that are denoted $\widehat{\Delta \left[ X\right]}_{G_j^N}$, $\widehat{\Delta \left[ Y\right]}_{G_j^N}$, $\widehat{\Delta \left[ X,Y\right]}_{G_j^N}$, $\widehat{\Delta \left[ X\right]}_{(I_Y)_j^N}$ and $\widehat{\Delta \left[ Y\right]}_{(I_X)_j^N}$ in the following. The underlain idea is to approximate the continuous random processes $(\sigma_t^X\sigma_t^Y\rho_t)^2$, $(\sigma_t^X)^2$ and $(\sigma_t^Y)^2$, or rather their time-transformed versions, by locally constant functions. This construction leads to time-adjusted histogram estimators  
\begin{subequations}
\begin{align}\label{avarest1}
\hat I_1&=\sum_{j=1}^{K_N}\left(\frac{\widehat{\Delta\left[ X,Y\right]}_{G_j^N}}{\Delta G_j^N}\right)^2\frac{G^N(T)}{K_N}~&\mbox{for}~~~~~&\int_0^T G^{\prime}(t)(\sigma_t^X\sigma_t^Y\rho_t)^2\,dt,~~~~\\
\label{avarest2}
\hat I_2&=\sum_{j=1}^{K_N}\left(\frac{\widehat{\Delta\left[ X\right]}_{G_j^N}\widehat{\Delta\left[ Y\right]}_{G_j^N}}{\left(\Delta G_j^N\right)^2}\right)\frac{G^N(T)}{K_N}~&\mbox{for}~~~~~&\int_0^T G^{\prime}(t)(\sigma_t^X\sigma_t^Y)^2\,dt,\\
\label{avarest3}
\hat I_3&=\sum_{j=1}^{K_N}\left(\frac{\widehat{\Delta\left[ X\right]}_{(I_Y)_j^N}}{\Delta (I_Y)_j^N}\right)\frac{I_Y^N(T)}{K_N}~&\mbox{for}~~~~~&\int_0^TI_Y^{\prime}(t)(\sigma_t^X)^2\,dt,\\
\label{avarest4}
\hat I_4&=\sum_{j=1}^{K_N}\left(\frac{\widehat{\Delta\left[ Y\right]}_{(I_X)_j^N}}{\Delta (I_X)_j^N}\right)\frac{I_X^N(T)}{K_N}~&\mbox{for}~~~~~&\int_0^TI_X^{\prime}(t)(\sigma_t^Y)^2\,dt.\end{align}
\end{subequations}
\begin{prop}\label{avarest}
The asymptotic variances \eqref{avarmulti} and \eqref{avarone} of the generalized multiscale estimator \eqref{multiscale} and the one-scale subsampling estimator \eqref{one-scale} with $M_N=c_{multi}N^{\nicefrac{1}{2}}$ and $i_N=c_{sub}N^{\nicefrac{2}{3}}$, can be estimated consistently by
\begin{subequations}
\begin{align}\label{avarestmulti}
\widehat{\AVAR}_{multi}&=\notag\left(c_{multi}^{-3}\left(24+12\,\frac{I_X^N(T)+I_Y^N(T)}{T}\right)+\frac{12}{5}c_{multi}^{-1}\right)\widehat{\eta_X^2}\widehat{\eta_Y^2}\\ 
&+c_{multi}\frac{26}{35}T\left(\hat I_1+\hat I_2\right)+c_{multi}^{-1}\frac{12}{5}\left(\widehat{\eta_Y^2}(1+\hat I_3)+\widehat{\eta_X^2}(1+\hat I_4)\right)~,
\end{align}
\begin{align}\label{avarestsub}
\widehat{\AVAR}_{sub}&=c_{sub}^{-2}4\widehat{\eta_X^2}\widehat{\eta_Y^2}
+c_{sub}\frac{2}{3}\left(\hat I_1+\hat I_2\right)~,
\end{align}
\end{subequations}
where $\hat I_1$-$\hat I_4$ are the estimators \eqref{avarest1}-\eqref{avarest4} and
\begin{align}\label{noisevarest}
\widehat{\eta_X^2}=(2n)^{-1}\sum_{i=1}^n(\Delta X_{t_i})^2~,~\widehat{\eta_Y^2}=(2m)^{-1}\sum_{j=1}^m(\Delta Y_{\tau_j})^2~.
\end{align}
\end{prop}
\begin{remark}
Convergence rates of the estimators \eqref{avarestmulti} and \eqref{avarestsub} for the asymptotic variances depend on the smoothness of $\sigma^X,\sigma^Y$ and $\rho$. For current stochastic volatility models as the Heston model, they are $N^{\nicefrac{1}{5}}$-consistent when choosing $K_N=c_KN^{\nicefrac{1}{5}}$ for a constant $c_K$ and $M_N\sim N^{\nicefrac{3}{5}}$ for the binwise multiscale estimators.\\
In the absence of noise, a consistent estimator for the asymptotic variances $2T\int_0^TG_X^{\prime}(t)(\sigma_t^X)^4dt$ of the realized volatility  has been proposed in \cite{bn3} as $(2n/3)\sum_{i=1}^n(\Delta X_{t_i})^4$. In the bivariate synchronous setting $(n/2)\sum_{i=1}^{n-1}(\Delta X_{t_i})^2\left((\Delta Y_{t_i})^2+(\Delta Y_{t_{i+1}})^2\right)$ is a convenient estimator. Consistency can be proved with It\^{o}'s formula and partial integration and comprehended by the analogy to a bivariate Gaussian distribution $(X,Y)\sim\mathbf{N}(0,\Sigma)$ with a covariance matrix $\Sigma$ with entries $\sigma_X^2,\sigma_Y^2,\rho\sigma_X\sigma_Y$. Then, $\E X^4=3\sigma_X^4$ and $\E\left[X^2Y^2\right]=2\rho^2\sigma_X^2\sigma_Y^2+\sigma_X^2\sigma_Y^2$ hold true.\\
In the noisy case smoothed versions of the estimators (using multiscale or alternative methods) are adequate (cf.\,\cite{kinnepoldi}).
However, in the non-synchronous non-noisy setting, there is no direct extension available and for that reason we have made the effort to construct the consistent histogram-based estimators \eqref{avarest1}-\eqref{avarest4} above.
\end{remark}
\section{Discussion and application\label{sec:6}}
\subsection{A case study}
We have learned by now that in a synchronous setting the special version of the central limit Theorem \ref{cltmulti} from Proposition \ref{cltsynchronous} holds true. Since asymptotics of the estimator \eqref{multiscale} not hinge on interpolations in the signal term the same central limit theorem applies in the case of intermeshed sampling introduced in Section \ref{sec:3}. Now we focus on observation schemes that arise as realizations of two homogeneous Poisson processes that are mutually independent and independent of the processes $\tilde X$ and $\tilde Y$. Although this model can be criticized for its flaw that sampling schemes of two correlated processes follow two independent processes and time homogeneity, what might seem to be  unrealistic in financial applications, independent and homogeneous Poisson sampling constitutes the most commonly used model in this research area (cf.\,\cite{zhang}, \cite{hy} among others) and appertains to show that the general form of \eqref{avarmulti} is tractable.\\
Let $\tilde n^{(n)}(t)$ and $\tilde m^{(n)}(t)$ be sequences of two independent homogeneous Poisson processes with parameters $Tn/\theta_1$ and $Tn/\theta_2$ ($n\in\N$), such that the waiting times between jumps of $\tilde n^{(n)}$ and $\tilde m^{(n)}$ are exponentially distributed with expectations $\E\left[\Delta t_i^{(n)}\right]=\theta_1/n$ and $\E\left[\Delta \tau_j^{(n)}\right]=\theta_2/n~,i\in\N,j\in\N$. In this case 
$$\Delta T_k^{(n)}\sim F(t)=1-\exp{\left(-\frac{tn}{\theta_1}\right)}-\exp{\left(-\frac{tn}{\theta_2}\right)}+\exp{\left(-tn\left(\frac{1}{\theta_1}+\frac{1}{\theta_2}\right)\right)}~,k\in\N\,,$$
\begin{align}\label{pp4}
I_X^N(t)\stackrel{p}{\longrightarrow} \frac{\theta_1\theta_2t}{(\theta_1+\theta_2)^2}~,~I_Y^N(t)\stackrel{p}{\longrightarrow} \frac{\theta_1\theta_2t}{(\theta_1+\theta_2)^2}~~\left(=\frac{z\,t}{(z+1)^2}~~\text{if}~~\theta_1=z\theta_2\right)~,\end{align}
hold true and we derive the following Poisson sampling version of Theorem \ref{cltmulti}:
\begin{cor}\label{cltmultipoiss}
On the Assumptions \ref{eff} and \ref{e}, the generalized multiscale estimator \eqref{multiscale} with noise-optimal weights \eqref{optimalweights}, and $M_N=c_{multi}\cdot \sqrt{N}$, converges conditionally on the independent Poisson sampling scheme with $0<\theta_1<\infty$ and $0<\theta_2<\infty$ stably in law with rate $N^{\nicefrac{1}{4}}$ to a mixed normal limit:
\begin{align*}N^{\nicefrac{1}{4}}\left(\widehat{\left[ X,Y\right]}_T^{multi}-\left[ X,Y\right]_T\right)\stackrel{st}{\rightsquigarrow}\mathbf{N}\left(0,\AVAR_{multi}^{poiss}\right)\end{align*}
with the asymptotic variance
\begin{align}\notag \AVAR_{multi}^{poiss}&=c_{multi}^{-3}\,\left(24+12\,\frac{2\theta_1\theta_2}{(\theta_1+\theta_2)^2}\right)\eta_X^2\eta_Y^2+c_{multi}^{-1}\,\frac{12\eta_X^2\eta_Y^2}{5}\\ 
&\; +c_{multi}\,\frac{26}{35}\int_0^T2\left(1-\frac{2\theta_1^2\theta_2^2}{\theta_1^2\theta_2^2+(\theta_1^2+\theta_2^2)(\theta_1+\theta_2)^2}\right)(\sigma_t^X\sigma_t^Y)^2(1+\rho_t^2)\,dt\\
\notag &\; +c_{multi}^{-1}\frac{12}{5}\left(\eta_Y^2\int_0^T(1+\frac{\theta_1\theta_2}{\theta_1+\theta_2})(\sigma_t^X)^2\,dt\,+\eta_X^2\int_0^T(1+\frac{\theta_1\theta_2}{\theta_1+\theta_2})(\sigma_t^Y)^2\,dt\right)~.\end{align}
The asymptotic variance of the $N^{\nicefrac{1}{6}}$-consistent one-scale estimator becomes
\begin{align}\notag \AVAR_{sub}^{poiss}&=c_{sub}^{-2}\,4\eta_X^2\eta_Y^2\\
&\; +c_{sub}\,\frac{2}{3}\int_0^T2\left(1-\frac{2\theta_1^2\theta_2^2}{\theta_1^2\theta_2^2+(\theta_1^2+\theta_2^2)(\theta_1+\theta_2)^2}\right)(\sigma_t^X\sigma_t^Y)^2(1+\rho_t^2)\,dt
~.\end{align}
\end{cor}
The order for the supremum of time instants in Assumption \ref{grid} holds in probability and the proof in \ref{sec:7} stays valid. A Poisson sampling version of Theorem \ref{HYclt} is given in \cite{asyn}, where the stochastic limit of $G^N$ is deduced using the distribution of a maximum of two exponentials stated above.
\subsection{A bridge between the noisy and the non-noisy setup}
So far we considered noise variances not dependent on $N$, but from an applied point of view there is interest in the case where the noise level may vary with $N\sim n \sim m$. The primary motivation to accommodate dependence of the noise on the sample size in the model originates from the economic background. Empirical studies of (ultra) high-frequency financial data suggest to rather model the observed log-prices as sum of a latent semimartingale and noise for that the variance decreases in $N$ as reported in \cite{kalnina} and \cite{distaso}, among others.\\
If an estimation approach uses previous-tick interpolations, as the one proposed in \cite{bn1}, these methods are not accurate for that setting any more. This also becomes apparent in the simulation study in \cite{bibinger} when the performance of the estimators is compared for varying noise levels. The generalized multiscale estimator is not biased due to asynchronicity and passes over to the Hayashi-Yoshida estimator \eqref{HY} for $M_N=1$, a $\sqrt{N}$-consistent estimator in the complete absence of noise. For that reason, our estimation method achieves an improved convergence rate in the model with decreasing noise variances. The next Corollary is obtained by a direct extension of the proof of Theorem \ref{cltmulti} in \ref{sec:7} when replacing the moments of the noise processes. A similar extension for the one-scale estimator where we obtain the rate $N^{\frac{1}{6}+\frac{\alpha}{3}}$ for a subsample frequency $i_N=c_{sub}N^{\frac{2}{3}(1-\alpha)}$ holds analogously.
\begin{cor}
Consider the model imposed by Assumption \ref{eff}, \ref{grid} and \ref{e}, but with noise variances $\eta_X^2(N)=\zeta_XN^{-\alpha}\,,\,\eta_Y^2=\zeta_YN^{-\alpha}\,,\,0<\alpha<1$ and constants $0<\zeta_X<\infty\,,\,0<\zeta_Y<\infty$. The generalized multiscale estimator \eqref{multiscale} with $M_N=c_{multi}N^{\frac{1}{2}-\frac{\alpha}{2}}$ and optimal weights \eqref{optimalweights} converges stably in law to a mixed Gaussian limit:
\begin{align}N^{\frac{1}{4}+\frac{\alpha}{4}}\left(\widehat{\left[ X,Y\right]}_T^{multi}-\left[ X,Y\right]_T\right)\stackrel{st}{\rightsquigarrow}\mathbf{N}\left(0,\AVAR_{multi}^*\right)\end{align}
with the asymptotic variance
\begin{align*}\notag \AVAR_{multi}^*&=c_{multi}^{-3}\left(24+12\,\frac{I_X(T)+I_Y(T)}{T}\right)\zeta_X\zeta_Y+c_{multi}^{-1}\,\frac{12\zeta_X\zeta_Y}{5}\\ 
&\; +c_{multi}\frac{26}{35}T\int_0^TG^{\prime}(t)(\sigma_t^X\sigma_t^Y)^2(1+\rho_t^2)\,dt\\ &\; +c_{multi}^{-1}\frac{12}{5}\left(\zeta_Y\int_0^T(1+I^{\prime}_Y(t))(\sigma_t^X)^2\,dt\,+\zeta_X\int_0^T(1+I^{\prime}_X(t))(\sigma_t^Y)^2\,dt\right)~.\end{align*}
\end{cor}
Incorporating a pure previous-tick interpolation strategy as in \cite{bn1}, one cannot gain an improved convergence rate in that setting due to the bias by non-synchronicity effects.
\subsection{Application}
\subsubsection{Note on modeling assumptions}
For an application to financial time series data the conditions imposed by Assumptions \ref{eff}-\ref{e} seem to be restrictive and the model will not describe stylized facts of the data in an adequate way. In particular relaxing the i.\,i.\,d.\,assumption on the noise is important. On the other hand, the underlying model in the sections before is convenient to establish an asymptotic distribution theory and ascertains a closed-form expression for the asymptotic variance. Yet, the generalization for serially dependent observation errors as carried out in \cite{zhangmykland2} for the one-dimensional case is possible for the generalized multiscale estimator \eqref{multiscale} as well. On the assumption of stationary strong mixing noise processes the multiscale estimator remains consistent and rate-optimal without any further adjustment. The analysis for the synchronous case can be adopted from \cite{zhangmykland2}, but for the general non-synchronous noisy setting a closed-form expression of the asymptotic variance and a corresponding limit theorem is not available. Furthermore, the condition that noise processes are mutually independent can be relaxed if one wants to allow some correlation $\E\left[\epsilon_{t_i}^X\epsilon_{\tau_j}^Y\right]=\eta_{X,Y}^{i,j}$ for $t_i$ and $\tau_j$ located near each other. In any case the generalized multiscale estimator remains asymptotically unbiased and $N^{\nicefrac{1}{4}}$-consistent. This does not necessarily hold true for the one-scale estimator that would have to be bias-corrected as the TSRV estimator by \cite{zhangmykland} in the one-dimensional case.\\
More general efficient processes including jumps can be covered in the model when we combine the method to a two-stage approach as presented by \cite{fanwang} for the one-dimensional estimation approach. Considering semimartingales
$$X_t=\int_0^t\mu_s^X\,ds+\int_0^t\sigma_s^XdB_s^X+\sum_{l=1}^{J_t^X}L_l^X~,~Y_t=\int_0^t\mu_s^Y\,ds+\int_0^t\sigma_s^YdB_s^Y+\sum_{l=1}^{J_t^Y}L_l^Y~,$$
with locally bounded drifts, continuous volatilities and counting processes $J_t^X,J_t^Y$ counting the jumps of $X$ and $Y$ with jump sizes $L_l^X,l=1,\ldots,J_t^X$ and $L_l^Y,l=1,\ldots,J_t^Y$, respectively, the generalized multiscale estimator converges in probability to the total quadratic covariation 
$$\int_0^T\rho_t\sigma_t^X\sigma_t^Y+\sum_{0\le s\le t}\Delta X_s\Delta Y_s$$
where $\Delta X_s=X_s-X_{s,-},\Delta Y_s=Y_s-Y_{s,-}$, and the second addend is the sum of the simultaneous co-jumps. When one is interested in disentangling the continuous part from the jumps, one convincing possibility following \cite{fanwang} is to use wavelet methods that locate the jumps in the sample paths, estimate jump sizes and afterwards use our estimation approach for the validated observations.\\
In conclusion the generalized multiscale estimator \eqref{multiscale} is capable for usage in various applications. In \cite{bibinger}, we have approved that the estimator performs well and has satisfying finite sample size features in simulations including serial dependent noise and typical stochastic volatility models.
\subsubsection{Choice of tuning parameters}
An implementation of the estimation approach requires first a rule to choose tuning parameters. We provide an accurate algorithm to implement the estimators and to obtain estimates for their asymptotic variances.
\addtocounter{figure}{-1} 
\renewcommand{\figurename}{Algorithm}
\begin{figure}[t]\fbox{
\begin{minipage}[c]{\textwidth} 
\renewcommand{\baselinestretch}{.75}\normalsize
\begin{itemize}\item Choose a priori $L$ and calculate pilot estimator $\widehat{\AVAR}_{multi}^{p}$ with 
\small{
\begin{align*}\widehat{\AVAR}_{dis}^{p}\hspace*{-.1cm}=\hspace*{-.1cm}\frac{N}{2}\hspace*{-.1cm}\sum_{k=1}^{\lfloor N/L\rfloor}\left(\left(\tilde X_{g_{kL}}-\tilde X_{l_{(k-1)L+1}}\right)^2\hspace*{-.1cm}+\hspace*{-.1cm}\left(\tilde X_{g_{kL+2}}-\tilde X_{l_{(k-1)L+3}}\right)^2\right)\hspace*{-.1cm}\left(\tilde Y_{\gamma_{kL}}-\tilde Y_{\lambda_{(k-1)L+1}}\right)^2\end{align*} }
\normalsize
and $\widehat{\left[ X\right]}_T^{p}=\sum_{j=kL,k\ge 1}\left(\tilde X_{t_j}-\tilde X_{t_{j-L}}\right)^2$
and $\widehat{\left[ Y\right]}_T^{p}$ analogously and $I_X(t)\equiv I_X^N(T)$ and $I_Y(t)\equiv I_Y^N(T)$. Calculate $\widehat{\eta_X^2}$ and $\widehat{\eta_Y^2}$ according to \eqref{noisevarest}.
\item Use pilot estimates to estimate optimal constant(s) 
\small{\begin{align}\label{cest}\hat c_{multi}^{(p)}=\left(\frac{-\widehat{\AVAR}_{cross,\text{n}}^{p}+\sqrt{\left(\widehat{\AVAR}_{cross,\text{n}}^{p}\right)^2+12\widehat{\AVAR}_{dis}^{p}\widehat{\AVAR}_{\text{n}}^{p}}}{6\widehat{\AVAR}_{\text{n}}^{p}}\right)^{-\nicefrac{1}{2}}\end{align}}
\normalsize
and $\hat c_{sub}^{(p)}=\sqrt[3]{2\widehat{\AVAR}_{\text{n},sub}^{p}/\widehat{\AVAR}_{dis,sub}^{p}}$.
\item Calculate $\hat I_1-\hat I_4$, given in \eqref{avarest1}-\eqref{avarest4}, with 
$$K_N=\sqrt{\hat c_{multi}^{(p)}}N^{\nicefrac{1}{5}}~~\text{bins and}~~M_N(j)=\left(\hat c_{multi}^{(p)}\right)^{\nicefrac{5}{4}}N^{\nicefrac{3}{5}}~\forall\,j~.$$
\item Estimate asymptotic variance with $\hat I_1-\hat I_4$ and $\widehat{\eta_X^2},\widehat{\eta_Y^2}$ and determine $\hat c_{multi}$ and $\hat c_{sub}$ with the above given formulae.
\item Calculate the generalized multiscale estimator \eqref{multiscale} with optimal weights \eqref{optimalweights} (and the one-scale subsample estimator)
with $M_N=\hat c_{multi}\sqrt{N}$ (and $i_N=\hat c_{sub}N^{\nicefrac{2}{3}}$)
\end{itemize}
\end{minipage}}
\caption{\label{AAVAR}Algorithm for the estimation procedure.}
\end{figure}
\renewcommand{\figurename}{Figure}
\renewcommand{\baselinestretch}{1}\normalsize\noindent
One plausible selection of the constants $c_{multi}=\sqrt{N}/M_N$ and $c_{sub}=N^{\nicefrac{2}{3}}/i_N$ can be derived as solutions of the minimization problems of the asymptotic variances. This leads to formula \eqref{cest} in Algorithm \ref{AAVAR}. \\
The tactic of Algorithm \ref{AAVAR} is the following:  Evaluate a pilot estimate $\hat c_{multi}^{(p)}$ for $c_{multi}$ as solution of formula \eqref{cest} inserting a sparse-sampled estimator for the signal term. Then set up the estimation of the asymptotic variances involving the estimators \eqref{avarest1}-\eqref{avarest4}. Take $M_N^b=c_{multi}^b\sqrt{NK_N}$ fixed for the multiscale estimators on all bins and set $M_N^b=c_{multi}^b\sqrt{Nc_{K}N^{\nicefrac{1}{5}}}$ where $c_K=K_N^{-1}N^{\nicefrac{1}{5}}$. This selection is optimal for common volatility models. We obtain $c_{multi}^bc_K^{-\nicefrac{1}{2}}=\hat c_{multi}^{(p)}$ and from the orders of the different errors of the histogram estimators $c_{multi}^b=c_K^{\nicefrac{5}{2}}$.
Hence, $c_{multi}^b=\left(c_{multi}^{(p)}\right)^{\nicefrac{5}{4}}$ and $c_K=\sqrt{c_{multi}^{(p)}}$ is derived. Using estimators \eqref{avarest1}-\eqref{avarest4}, we calculate estimates for the addends of the asymptotic variance and $\hat c_{multi}$ according to formula \eqref{cest} again and $M_N=\lceil \hat c_{multi}\sqrt{N}\, \rceil$ is used for the final estimator. It turns out that this strategy it quite robust to the a priori chosen sparse-sample frequency that can be chosen under the impression of usual diagnostic tools as signature plots and acfs. 
\subsubsection{Simulation and data analysis check}
As completion to the detailed simulation study in the supplementary material to \cite{bibinger} we investigate the performance of Algorithm \ref{AAVAR} and the histogram-estimators \eqref{avarest1}-\eqref{avarest4} here. For this purpose we simulate from an simple Brownian motion model with zero drifts and constant volatilities $\sigma^X=\sigma^Y=1$ and $\rho=1/2$ with equal noise variances $\eta^2$. Sampling schemes are generated by independent time-homogeneous Poisson sampling with 30.000 expected observation for both processes on $[0,1]$. Results of the estimates are listed in Table \ref{stab2}.
\begin{table}
\caption{\label{stab2}Estimators \eqref{avarest1}-\eqref{avarest4}, estimators for the asymptotic variances of the multiscale \eqref{avarmulti} and the one-scale estimator \eqref{avarone}, calculated asymptotic variances and estimates for the quadratic covariation. The estimates are given $\pm$ empirical standard deviations.}
\begin{center}\hspace*{-0.125cm}
\begin{tabular}{|l|ccccc|}
\hline noise var.\,$\eta^2$&0.0001&$(0.001/\sqrt{10})$&0.001&$(0.01/\sqrt{10})$&0.01\\
\hline
$\hat I_1$ \hspace*{-0.25cm}&\hspace*{-0.15cm}$0.392\pm 0.038$ &$0.390\pm 0.047$ &$0.394\pm 0.073$&$0.413\pm 0.144$&$0.423\pm 0.128$\\
$\hat I_2$ \hspace*{-0.25cm}&\hspace*{-0.15cm}$1.557\pm 0.067$ &$1.552\pm 0.085$ &$1.538\pm 0.141$&$1.529\pm 0.276$&$1.462\pm 0.230$\\
$\hat I_3$ \hspace*{-0.25cm}&\hspace*{-0.15cm}$0.250\pm 0.007$ &$0.249\pm 0.010$ &$0.249\pm 0.016$&$0.247\pm 0.031$&$0.234\pm 0.085$\\
$\hat I_4$ \hspace*{-0.25cm}&\hspace*{-0.15cm}$0.250\pm 0.007$ &$0.249\pm 0.009$ &$0.249\pm 0.016$&$0.247\pm 0.030$&$0.233\pm 0.083$\\
$\widehat{\AVAR}_{multi}$ \hspace*{-0.25cm}&\hspace*{-0.15cm}$0.090\pm 0.003$ &$0.143\pm 0.005$ &$0.246\pm 0.015$&$0.434\pm 0.050$&$0.778\pm 0.082$\\
$\AVAR_{multi}$ \hspace*{-0.25cm}&\hspace*{-0.15cm}$0.0663$ &$0.1185$ &$0.2159$&$0.3774$&$0.6737$\\
$\widehat{\AVAR}_{sub}$ \hspace*{-0.25cm}&\hspace*{-0.15cm}$0.0086\pm 0.0002$ &$0.017\pm 0.001$ &$0.037\pm 0.002$&$0.080\pm 0.009$&$0.157\pm 0.017$\\
$\AVAR_{sub}$ \hspace*{-0.25cm}&\hspace*{-0.15cm}$0.0077$ &$0.0166$ &$0.0357$&$0.0768$&$0.1656$\\
$\widehat{\left[ X,Y\right]}_T^{multi}$ \hspace*{-0.25cm}&\hspace*{-0.15cm}$0.501\pm 0.024$ &$0.499\pm 0.029$ &$0.498\pm 0.038$&$0.499\pm 0.049$&$0.501\pm 0.065$\\
$\widehat{\left[ X,Y\right]}_T^{sub}$\hspace*{-0.25cm} &\hspace*{-0.15cm}$0.500\pm 0.022$ &$0.499\pm 0.028$ &$0.499\pm 0.042$&$0.500\pm 0.058$&$0.503\pm 0.074$\\
\hline
\end{tabular}
\end{center}
\end{table}
In addition, our approach is tested in an application to EUREX future tick-data taken from a database provided by the Research Data Center (RDC) of the CRC 649 \lq Economic Risk\rq~in Berlin. We aim at estimating integrated covolatilities between the four financial securities with the highest tick-frequencies in the database. These are the Euro-Bund Future (FGBL), that is based on a notional long-duration debt instrument issued by the Federal Republic of Germany, the Euro-Bobl Future (FGBM), a likewise medium-duration contract, and futures on the EURO STOXX 50 (FESX) and the German DAX (FDAX).\\
We apply the procedure with Algorithm \ref{A} and Algorithm \ref{AAVAR} to a suitable filtered dataset and give the results for two days in Table \ref{tabdat1} with associated estimated optimal multiscale frequencies. The estimates for the quadratic covariations are given $\pm$ estimated standard deviation from the estimated asymptotic variances and the pertaining $N$ for each pair. Although there are characteristics of the dataset not in accordance with the model assumptions, above all price discreteness and the fact that most returns are zero, the estimation approach passes this intuition check. Since the ESX and the DAX share 13 companies constituting c.\,28.5\% weighting in the ESX and c.\,72.4\% in the DAX there is a big systematic positive correlation between both and we presume that there is as well a high correlation between the two debt instruments which is both revealed by the estimates. On 09/11/2001, there has been an tremendous impact so that FGBL/FGBM have increased and the FESX/FDAX decreased. For that day we have significantly negative integrated covolatilities between debt instruments and stock indices which is not the case for the common trading day in comparison. As answer to the great amount of zero returns the only adjustment of the method that we undertake is to estimate noise variances in \eqref{noisevarest} by dividing the realized volatility by twice the number of non-zero returns instead of all returns.\\
To sum up, the estimation approach based on Algorithm \ref{A}, a multiscale extension of subsampling and Algorithm \ref{AAVAR} provides a convincing method to obtain integrated covolatility estimates for very general high-frequency data.
\begin{table}[h]
\begin{center}
\begin{tabular}{|l||c|c|c|c|}
\hline
 $\widehat{\left[ X,Y\right]}_T^{multi}$~(MSFR) & FGBL & FGBM & FESX & FDAX \\
 \hline\hline 
 FGBL & $5.46\pm .87$ ~(7)& $3.32\pm .74$~(7) & $0.98\pm .99$~(9)& $0.52\pm .74$~(6)\\
 FGBM & \cellcolor{gray} & $2.42\pm .47$~(6) & $0.78\pm .69$ ~(7)& $0.64\pm .51$~(4)\\
 FESX & \cellcolor{gray} & \cellcolor{gray} &$68.26\pm 7.1$~(10)&$29.39\pm 2.4$~(7)\\
 FDAX & \cellcolor{gray} & \cellcolor{gray} & \cellcolor{gray} &$61.43\pm 3.5$~(5)\\
 \hline
\end{tabular}\\
\begin{tabular}{|l||c|c|c|c|}
\hline
 $\widehat{\left[ X,Y\right]}_T^{multi}$~(MSFR) & FGBL & FGBM & FESX & FDAX \\
 \hline\hline 
 FGBL & $27.89\pm 4.1$ ~(15)& $12.94\pm 1.7$~(12) & $-52.55\pm 18$~(8)& $-34.13\pm 11$~(7)\\
 FGBM & \cellcolor{gray} & $18.10\pm 2.5$~(8) & $-26.44\pm 14$ ~(6)& $-25.01\pm 28$~(3)\\
 FESX & \cellcolor{gray} & \cellcolor{gray} &$3070\pm 172$~(6)&$757\pm 15$~(4)\\
 FDAX & \cellcolor{gray} & \cellcolor{gray} & \cellcolor{gray} &$1870\pm 94$~(4)\\
 \hline
\end{tabular}
\end{center}
\caption{\label{tabdat1}Estimates for integrated covolatilities, standard deviations ($\cdot$ $10^6$) and used multiscale frequencies for 10/01/2008 (top) and 09/11/2001 (bottom).}
\end{table}



\appendix
 \section{Proof of Theorem \ref{cltmulti}\label{sec:7}}
\subsection{Error due to noise and choosing the weights\label{subsec:4.3.1}}
\renewcommand{\appendixname}{}
The error due to microstructure noise of the generalized multiscale estimator is given by \begin{align*}&\sum_{i=1}^{M_N}\frac{\alpha_{i,M_N}}{i}\sum_{j=i}^N\left(\epsilon_{g_j}^X-\epsilon_{l_{j-i+1}}^X\right)\left(\epsilon_{\gamma_j}^Y-\epsilon_{\lambda_{j-i+1}}^Y\right)
=\sum_{i=1}^{M_N}\frac{\alpha_{i,M_N}}{i}\Big(\sum_{j=1}^N\left(\epsilon_{g_j}^X\epsilon_{\gamma_j}^Y+\epsilon_{l_j}^X\epsilon_{\lambda_{j}}^Y\right)\\ &\hspace*{3.795cm}-\sum_{j=i}^N\left(\epsilon_{g_j}^X\epsilon_{\lambda_{j-i+1}}^Y+\epsilon_{\gamma_j}^Y\epsilon_{l_{j-i+1}}^X\right)-\sum_{j=1}^{i-1}\epsilon_{g_j}^X\epsilon_{\gamma_j}^Y-\sum_{j=N-i+1}^{N}\epsilon_{l_j}^X\epsilon_{\lambda_j}^Y\Big).\end{align*}
Additionally to the standardization condition
\begin{equation}\label{c1}\sum_{i=1}^{M_N}\alpha_{i,M_N}=1~,\tag{C1}\end{equation}
that is necessary for asymptotic unbiasedness and consistency,
we now impose the auxiliary condition
\begin{equation}\label{c2}\sum_{i=1}^{M_N}\frac{\alpha_{i,M_N}}{i}=0~,\tag{C2}\end{equation}
on the weights which assures that the leading term of the noise error equals zero. Hence, there remain three uncorrelated addends in the error induced by microstructure noise.

\begin{prop}\label{propnoise1}
Let Assumptions \ref{grid} and \ref{grid4} on the observation times and Assumption \ref{e} on the observation errors hold true.
The asymptotic variance of the term
$$-\sum_{i=1}^{M_N}\frac{\alpha_{i,M_N}}{i}\sum_{j=i}^N\left(\epsilon_{g_j}^X\epsilon_{\lambda_{j-i+1}}^Y+\epsilon_{\gamma_j}^Y\epsilon_{l_{j-i+1}}^X\right)$$
is minimized by the weights
\begin{align}\label{optimalweights}
\alpha_{i,M_N}^{opt}=\left(\frac{12i^2}{(M_N^3-M_N)}-\frac{6i}{(M_N^2-1)}-\frac{6i}{(M_N^3-M_N)}\right)=\frac{12i^2}{M_N^3}-\frac{6i}{M_N^2}\left(1+\KLEINO(1)\right)\end{align}
as $M_N,N\rightarrow\infty$ and $M_N/N \rightarrow 0$ with $N=\KLEINO\left(M_N^4\right)$.
The following asymptotic normality result holds true:
\begin{align}\sqrt{\frac{M_N^3}{N}}\left(\sum_{i=1}^{M_N}\frac{\alpha_{i,M_N}^{opt}}{i}\sum_{j=i}^N\left(\epsilon_{g_j}^X\epsilon_{\lambda_{j-i+1}}^Y-\epsilon_{\gamma_j}^Y\epsilon_{l_{j-i+1}}^X\right)\right)\rightsquigarrow\mathbf{N}\left(0\,,\,\AVAR_{\text{noise}}\right)~,\end{align}
with the asymptotic variance
\begin{equation}\label{avarnoise1}\AVAR_{\text{noise}}=(24+12 (I_X(T)+I_Y(T)))\eta_X^2\eta_Y^2\end{equation}
with the functions $I_X$ and $I_Y$ defined in Assumption \ref{grid4}. The weak convergence also holds true conditionally given the paths of the efficient processes.
\end{prop}
\begin{proof}The term is centred and we illustrate it in the way
\begin{align*}-\sum_{i=1}^{M_N}\frac{\alpha_{i,M_N}}{i}\sum_{j=i}^N\left(\epsilon_{g_j}^X\epsilon_{\lambda_{j-i+1}}^Y+\epsilon_{\gamma_j}^Y\epsilon_{l_{j-i+1}}^X\right)
=-\sum_{j=1}^N\sum_{i=1}^{M_N\wedge j}\frac{\alpha_{i,M_N}}{i}\Big(\epsilon_{g_j}^X\epsilon_{\lambda_{j-i+1}}^Y\hspace*{2.5cm}\\\hspace*{2.5cm}(\1_{\{g_j=g_{j+1}\}}+\1_{\{g_j\ne g_{j+1}\}})
+\epsilon_{\gamma_j}^Y\epsilon_{l_{j-i+1}}^X(\1_{\{\gamma_j=\gamma_{j+1}\}}+\1_{\{\gamma_j\ne \gamma_{j+1}\}})\Big)~.\end{align*}
For fixed $i$ the addends of the inner sum are uncorrelated because $l_i\ne l_j$ and $\lambda_i\ne \lambda_j$ for all $i\ne j$. Consecutive right-end points $g_i,\gamma_i$ can be the same observation times instead, so that the inner sums are 2-dependent random variables.
Thus, the variance is given by
\begin{align*}&\var\left(\sum_{j=1}^{N}\frac{\alpha_{i,M_N}}{i}\sum_{i=1}^{M_N\wedge j}\left(\epsilon_{g_j}^X\epsilon_{\lambda_{j-i+1}}^Y-\epsilon_{\gamma_j}^Y\epsilon_{l_{j-i+1}}^X\right)\right)\\
&=\sum_{j=1}^N\sum_{i=1}^{M_N\wedge j}\left(\frac{\alpha_{i,M_N}}{i}\right)^22\eta_X^2\eta_Y^2+\sum_{j=1}^N\sum_{i=1}^{(M_N-1)\wedge (j-1)}\frac{\alpha_{i,M_N}\alpha_{i+1,M_N}}{i(i+1)}\eta_X^2\eta_Y^2(1-\1_{\{g_j\ne g_{j+1}\,,\,\gamma_j\ne\gamma_{j+1}\}})~.\end{align*}
The weights that minimize the first addend of the above variance and also the total variance asymptotically have been determined in \cite{bibinger}. Those weights are in line with the standard weights from \cite{zhang} in the univariate setting and correspond to a cubic kernel for the kernel estimator by \cite{bn1}.\\
Inserting the noise-optimal weights \eqref{optimalweights}, we can apply a central limit theorem for strong mixing triangular arrays from \cite{utev} to
$$\sqrt{\frac{M_N^3}{N}}\sum_{j=1}^{N}\frac{\alpha_{i,M_N}^{opt}}{i}\sum_{i=1}^{M_N\wedge j}\left(\epsilon_{g_j}^X\epsilon_{\lambda_{j-i+1}}^Y+\epsilon_{l_{j-i+1}}^X\epsilon_{\gamma_j}^Y\right)~.$$
The sequence of variances with the chosen weights according to \eqref{optimalweights}
\begin{align*}&\var\left(\sqrt{\frac{M_N^3}{N}}\sum_{j=1}^{N}\frac{\alpha_{i,M_N}^{opt}}{i}\sum_{i=1}^{M_N\wedge j}\left(\epsilon_{g_j}^X\epsilon_{\lambda_{j-i+1}}^Y+\epsilon_{l_{j-i+1}}^X\epsilon_{\gamma_j}^Y\right)\right)\\
&~=\frac{M_N^3}{N}\sum_{j=0}^N\sum_{i=1}^{M_N\wedge(j+1)}\left(\frac{\alpha_{i,M_N}^{opt}}{i}\right)^2\eta_X^2\eta_Y^2(3-\1_{\{g_j\neg_{j+1}\,,\,\gamma_j\ne\gamma_{j+1}\}})+\KLEINO(1)\\
&~\longrightarrow36\eta_X^2\eta_Y^2-12(1-(I_X(T)+I_Y(T)))\eta_X^2\eta_Y^2\end{align*}
converges to $\AVAR_{\text{noise}}$ on the Assumption \ref{grid4}.\\
Since the inner sums are 2-dependent and hence in particular $\phi$-mixing, the Lyapunov condition that holds obviously suffices to apply the central limit theorem of \cite{utev}. This completes the proof of the proposition.\end{proof}
Next, we consider the remaining addends of the error induced by microstructure frictions 
and insert the weights \eqref{optimalweights}:
\begin{prop}\label{propnoise2}On the Assumptions \ref{grid}, \ref{grid4} and \ref{e}, the following weak convergence to a centred normal distribution holds true:
\begin{align}\sqrt{M_N}\left(\sum_{i=1}^{M_N}\frac{\alpha_{i,M_N}^{opt}}{i}\left(\sum_{j=1}^{i-1}\epsilon_{g_j}^X\epsilon_{\gamma_j}^Y+\sum_{j=N-i+1}^N\epsilon_{l_j}^X\epsilon_{\lambda_j}^Y\right)\right)\rightsquigarrow\mathbf{N}\left(0,\frac{12}{5}\eta_X^2\eta_Y^2\right)~.\end{align}
This convergence also holds conditionally on the paths of the efficient processes.
\end{prop}
\begin{proof}
\begin{align*}\sum_{i=1}^{M_N}\frac{\alpha_{i,M_N}^{opt}}{i}\left(\sum_{j=1}^{i-1}\epsilon_{g_j}^X\epsilon_{\gamma_j}^Y+\sum_{j=N-i+1}^N\epsilon_{l_j}^X\epsilon_{\lambda_j}^Y\right)=\sum_{j=1}^{M_N-1}\left(\epsilon_{g_j}^X\epsilon_{\gamma_j}^Y+\epsilon_{l_{N-j}}^X\epsilon_{\lambda_{N-j}}^Y\right)\sum_{i=j+1}^{M_N}\frac{\alpha_{i,M_N}^{opt}}{i}
\end{align*}
Both addends are uncorrelated and treated analogously. We restrict ourselves to the proof for the first term. $\sqrt{M_N}\sum_{j=1}^{M_N-1}\epsilon_{g_j}^X\epsilon_{\gamma_j}^Y\sum_{i=j+1}^{M_N}\alpha_{i,M_N}^{opt}/i$ is the endpoint of a discrete centred martingale with respect to the filtration $\overline{\mathcal{A}}_{j}^{N}\:=\sigma\left(\epsilon_{t_k}^X|t_k\le g_j\,,\,X_{t_k}|0\le k\le n\right)\vee\sigma\left(\epsilon_{\tau_k}^Y|\tau_k\le \gamma_j\,,\,Y_{\tau_k}|0\le k\le m\right)$. Namely, since $g_j=g_{j-1}\Rightarrow \gamma_j> \gamma_{j-1}$ and $\gamma_j=\gamma_{j-1}\Rightarrow g_j> g_{j-1}$ analogously:
\begin{align*}
&\E\left[\sqrt{M_N}\epsilon_{g_l}^X\epsilon_{\gamma_l}^Y\sum_{i=l+1}^{M_N}\frac{\alpha_{i,M_N}^{opt}}{i}\Big|\overline{\mathcal{A}}_{l-1}^{N}\right]\\
&=\hspace*{-0.05cm}\sqrt{M_N}\hspace*{-0.05cm}\left(\1_{\{g_l=g_{l-1}\}}\E\left[\epsilon_{\gamma_l}^Y\right]\epsilon_{g_{l-1}}^X\hspace*{-0.05cm}+\hspace*{-0.05cm}\1_{\{\gamma_l=\gamma_{l-1}\}}\E\left[\epsilon_{g_l}^X\right]\epsilon_{\gamma_{l-1}}^Y\hspace*{-0.05cm}+\hspace*{-0.05cm}\1_{\{g_l\ne g_{l-1}\,,\,\gamma_l\ne\gamma_{l-1}\}}\E\left[\epsilon_{g_l}^X\right]\E\left[\epsilon_{\gamma_l}^Y\right]\right)\hspace*{-0.05cm}=\hspace*{-0.05cm}0\,.\end{align*}
A central limit theorem for martingale triangular arrays from \cite{hall} is applied, in particular the non-stable version of Corollary 3.\,1 (cf.\,the following remark in \cite{hall} and references cited therein).
The conditional Lindeberg condition can be verified by the stronger conditional Lyapunov condition. The proof of it is obtained by a similar calculation as the following one and we omit it. The conditional variance equals
\begin{align*}&M_N\sum_{j=1}^{M_N-1}\var\left(\epsilon_{g_j}^X\epsilon_{\gamma_j}^Y\sum_{i=j+1}^{M_N}\frac{\alpha_{i,M_N}^{opt}}{i}\Big|\mathcal{A}_{N,j-1}\right)=M_N\sum_{j=1}^{M_N-1}\left(\sum_{i=j+1}^{M_N}\frac{\alpha_{i,M_N}^{opt}}{i}\right)^2\\
&\hspace*{3cm}\times\Big(\eta_X^2\eta_Y^2\1_{\{g_j\ne g_{j-1}\,,\,\gamma_j\ne\gamma_{j-1}\}}+\left(\epsilon_{g_j}^X\right)^2\1_{\{g_j=g_{j-1}\}}\eta_Y^2+\left(\epsilon_{\gamma_j}^Y\right)^2\1_{\{\gamma_j=\gamma_{j-1}\}}\eta_X^2\Big)\\
&~\stackrel{p}{\longrightarrow}\frac{6}{5}\eta_X^2\eta_Y^2~.\end{align*}
We have used the formula $\sum_{j=1}^{M_N-1}\left(\sum_{i=j+1}{M_N}(\alpha_{i,M_N}^{opt}/i)\right)^2=(6/5)M_N^{-1}+\KLEINO(M_N^{-1})$.\\
\end{proof}
\subsection{Discretization errors of the estimators}
\begin{prop}\label{cltosdis}On the Assumptions \ref{eff}, \ref{grid} and \ref{grid3}, the discretization error of the one-scale subsampling estimator with subsampling frequency $i_N$ converges stably in law to a centred mixed normal limit as $i_N\rightarrow\infty, N\rightarrow\infty, i_N/N^{\alpha}\rightarrow 0$ for every $\alpha>2/3$:
\begin{align*}\sqrt{\frac{N}{i_N }}\left(\sum_{j=i}^N\left(X_{g_j}-X_{l_{j-i+1}}\right)\left(Y_{\gamma_{j}}-Y_{\lambda_{j-i+1}}\right)-\left[ X,Y\right]_T\right)\stackrel{st}{\rightsquigarrow}\mathbf{N}\left(0,\AVAR_{dis,sub}\right)~,\end{align*}
with asymptotic variance
\begin{align}\AVAR_{dis,sub}=\frac{2}{3}T\int_0^TG^{\prime}(t)(\sigma_t^X\sigma_t^Y)^2(1+\rho_t^2)\,dt~.\end{align}
\end{prop}
\begin{prop}\label{cltmsdis} On the Assumptions \ref{eff},\ref{grid} and \ref{grid3}, the discretization error of the generalized multiscale estimator with the noise-optimal weights given in \eqref{optimalweights} converges with rate $\sqrt{N/M_N}$ stably in law to a centred mixed Gaussian limit as $M_N\rightarrow\infty, N\rightarrow\infty, M_N/N^{\alpha}\rightarrow 0$ for every $\alpha>2/3$:
\begin{align*}\sqrt{\frac{N}{M_N}}\left(\sum_{i=1}^{M_N}\frac{\alpha_{i,M_N}^{opt}}{i}\sum_{j=i}^N\left(X_{g_j}-X_{l_{j-i+1}}\right)\left(Y_{\gamma_j}-Y_{\lambda_{j-i+1}}\right)-\left[ X,Y\right]_T\right)\stackrel{st}{\rightsquigarrow}\mathbf{N}\left(0,\AVAR_{dis,multi}\right)~,\end{align*}
with asymptotic variance
\begin{align}\label{avardis}\AVAR_{dis,multi}=\frac{26}{35}T\int_0^TG^{\prime}(t)(\sigma_t^X\sigma_t^Y)^2(1+\rho_t^2)\,dt~.\end{align}
\end{prop}
\subsubsection{Time-change in the asymptotic quadratic variation of time}
\begin{prop}\label{tcr}In the proof of a central limit theorem for the discretization error of the closest synchronous approximation $T_k^{(N)},\,k=0,\ldots,N$, of our generalized multiscale estimator \eqref{multiscale} on the Assumptions \ref{grid} and \ref{grid4}, we can additionally, without loss of further generality, assume that
\begin{align}\label{tcaqvt}\sum_{k=1}^N\left(\Delta T_k^{(N)}-\frac{T}{N}\right)^2=\KLEINO\left(N^{-1}\right)~.\end{align}
\end{prop}
\begin{remark}
From Assumptions \ref{eff} and \ref{grid}, we can deduce directly that the sum above is at most of order $N^{-1}$. The stronger assertion, that the closest synchronous approximation defined by the times $T_k^{(N)},\,k=0,\ldots,N$ introduced in paragraph \ref{subsec:HY} is close to equidistant sampling in the sense that the sum above is of smaller asymptotic order than $N^{-1}$, is derived by the concept of a time-change in the asymptotic quadratic variation of time from Assumption \ref{grid3}. For the proof of \eqref{tcaqvt} we refer to Lemma 1 from \cite{zhang} where this concept has been presented for the univariate multiscale approach and it directly carries over to the synchronous multivariate case.\\
On the Assumption \ref{grid3}, a transformation $g$ can be defined that maps the refresh times $T_k^{(N)}$ to values $g(T_k^{(N)})$, so that \eqref{tcaqvt} holds true for the transformed synchronous observation scheme. Thanks to the fact that the corresponding time-changed processes $L_{g(t)}$ and $M_{g(t)}$ fulfill Assumption \ref{eff} again and the transformed observation scheme Assumption \ref{grid}, we are able to prove a central limit theorem for the time-changed version of the discretization error if \eqref{tcaqvt} did not hold.\\
Since the resulting asymptotic variance will be invariant under the transformation $g$, the central limit theorem will analogously hold true for the original sampling scheme. Hence, no further restriction has to be made when assuming \eqref{tcaqvt}.
\end{remark}
\subsubsection{Discretization error of the closest synchronous approximation}
Note that it suffices to prove the foregoing limit theorems for the zero-drift case. Since our limit theorems are stable, asymptotic mixed normality is assured to hold for the general setting on Assumption \ref{eff}. Denote $L_t=\int_0^t\sigma_s^XdW_s^X$ and $M_t=\int_0^t\sigma_s^YdW_s^Y$ the continuous martingales that represent the efficient processes under the equivalent martingale measure after a Girsanov transformation. The Novikov condition has been imposed in Assumption \ref{eff} to allow for this transformation.\\
The asymptotic mixed normality result is implied as marginal distribution at $t=T$ of a limiting time-changed Brownian motion which is proven to be the stable weak limit of the process corresponding to the discretization error \eqref{diserr} with the theory of \cite{jacod1}. \\
We begin with the discretization error of the closest synchronous approximation of a one-scale subsampling estimator. 
\begin{prop}\label{cltosdissynpr}On the same assumptions as in Proposition \ref{cltosdis}, the continuous martingale 
\begin{align*}\mathfrak{D}_t^N\:=\sqrt{\frac{N}{i_N T}}&\left[\sum_{T_k\le t}\left(\Delta L_{T_k}+\int_{T_k}^t\sigma_s^XdW_s^X\right)\left(\sum_{l=1}^{i\wedge k}\left(1-\frac{l}{i}\right)\Delta M_{T_{k-l}}\right)\right. \\
&\left. +\sum_{T_k\le t}\left(\Delta M_{T_k}+\int_{T_k}^t\sigma_s^YdW_s^Y\right)\left(\sum_{l=1}^{i\wedge k} \left(1-\frac{l}{i}\right)\Delta L_{T_{k-l}}\right)\right]\end{align*}
for $t\in[0,T]$, where $\Delta\,\cdot\,_{T_k}=\,\cdot\,_{T_k}-\,\cdot\,_{T_{k-1}}$ is the backward difference operator, converges stably in law as $N\rightarrow\infty, i_N\rightarrow\infty, i_N/N\rightarrow 0$ to a limiting time-changed Brownian motion
\begin{align*}\mathfrak{D}^N_t\stackrel{st}{\rightsquigarrow}\int_0^t\sqrt{v_{\mathfrak{D}_s}}d{\mathfrak{W}}_s^{\bot}~,\end{align*}
where $\mathfrak{W}^{\bot}$ is independent of $\mathcal{F}$ and
\begin{align*} v_{\mathfrak{D}_s}=\frac{2}{3}  G^{\prime}(s)(\sigma_s^X\sigma_s^Y)^2(1+\rho_s^2)~.\end{align*}
\end{prop}
\noindent
\textit{Proof of Proposition \ref{cltosdissynpr}:}\\
The subscript of the subsampling frequency is omitted in the following proof and $C$ denotes a generic constant and $\delta_N=\sup_{i\in\{1,\ldots,N\}}{(T_{i}-T_{i-1})}$. \\
We apply a simplified martingale version of the stable central limit theorem 2--1 from \cite{jacod1}. For other applications and expositions of the theory from \cite{jacod1} we refer to \cite{poldilimit}, \cite{fukasawa} and \cite{asyn}. The above limit theorem is implied by the following three conditions:
\begin{subequations}
\begin{align}\label{qvjac}\left[\mathfrak{D}\right]_t\stackrel{p}{\longrightarrow}\int_0^tv_{\mathfrak{D}_s}\,ds~,\end{align}
\begin{align}\label{covjac}\left[\mathfrak{D},L\right]_t\stackrel{p}{\longrightarrow}0~,~\left[\mathfrak{D},M\right]_t\stackrel{p}{\longrightarrow}0,~\end{align}
\begin{align}\label{ortjac}\left[\mathfrak{D},L^{\bot}\right]_t\stackrel{p}{\longrightarrow}0~,~\left[\mathfrak{D},M^{\bot}\right]_t\stackrel{p}{\longrightarrow}0,~\end{align}
\end{subequations}
for all $t \in[0,T]$ and for all $M^{\bot}\in\mathcal{M}{\bot}$ and $L^{\bot}\in\mathcal{L}{\bot}$ that denote the set of $\left(\mathcal{F}\right)_t$-adapted bounded martingales orthogonal to $M$ and $L$, respectively.\\
Calculating the quadratic variation of $\mathfrak{D}_t^N$ yields
\begin{align*} \left[ \mathfrak{D}^N\right]_t&=\frac{N}{iT}\left[\sum_{T_k\le t}\left(\Delta\left[ L\right]_{T_k}\left(\sum_{l=1}^{i\wedge k}\left(1-\frac{l}{i}\right)\Delta M_{T_{k-l}}\right)^2+\Delta\left[ M\right]_{T_k}\left(\sum_{l=1}^{i\wedge k}\left(1-\frac{l}{i}\right)\Delta L_{T_{k-l}}\right)^2\right)\right.\\
 & \left. ~~~~~~~~~~+2\sum_{T_k\le t}\Delta\left[ L,M\right]_{T_k}\left(\sum_{l=1}^{i\wedge k}\left(1-\frac{l}{i}\right)\Delta L_{T_{k-l}}\right)\left(\sum_{l=1}^{i\wedge k}\left(1-\frac{l}{i}\right)\Delta M_{T_{k-l}}\right)\right]+\KLEINO_p(1) \displaybreak[0]\\
 &=~\frac{N}{iT}\left[\sum_{T_k\le t}\Delta\left[ L\right]_{T_k}\sum_{l=1}^{i\wedge k}\left(1-\frac{l}{i}\right)^2\left(\Delta M_{T_{k-l}}\right)^2+\sum_{T_k\le t}\Delta\left[ M\right]_{T_k}\right.
 \\& \left.~~~~~~~~\times \sum_{l=1}^{i\wedge k}\left(1-\frac{l}{i}\right)^2\left(\Delta L_{T_{k-l}}\right)^2
  +2\sum_{T_k\le t}\Delta\left[ L,M\right]_{T_k}\sum_{l=1}^{i\wedge k}\left(1-\frac{l}{i}\right)^2\left(\Delta L_{T_{k-l}}\right)^2\right]+\KLEINO_p(1) \displaybreak[0]\\
  &=~\frac{N}{iT}\left[\sum_{T_k\le t}\int_{T_{k-1}}^{T_k}(\sigma_s^X)^2ds\left(\sum_{l=1}^{i\wedge k}\left(1-\frac{l}{i}\right)^2\int_{T_{k-l-1}}^{T_{k-l}}(\sigma_s^Y)^2ds\right)\right.
 \\&~~~~~~~~~~~~~~~ \left.+\sum_{T_k\le t}\int_{T_{k-1}}^{T_k}(\sigma_s^Y)^2ds \left(\sum_{l=1}^{i\wedge k}\left(1-\frac{l}{i}\right)^2\int_{T_{k-l-1}}^{T_{k-l}}(\sigma_s^X)^2ds\right)\right.
 \\& \left.
 ~~~~~~~~~~~~~~~ +2\sum_{T_k\le t}\int_{T_{k-1}}^{T_k}\rho_s\sigma_s^X\sigma_s^Yds\left(\sum_{l=1}^{i\wedge k}\left(1-\frac{l}{i}\right)^2\int_{T_{k-l-1}}^{T_{k-l}}\rho_s\sigma_s^X\sigma_s^Yds\right)\right]+\KLEINO_p(1)\displaybreak[0]\\
  &{\underset{\text{\tiny{Lemma \ref{crucialh2}}}}{=}}~\frac{N}{iT}\sum_{T_k\le t}2(1+\rho_{T_{k-1}}^2)(\sigma_{T_{k-1}}^X\sigma_{T_{k-1}}^Y)^2\left(\Delta T_k\right)^2\sum_{l=1}^{i\wedge k}\left(1-\frac{l}{i}\right)^2+\KLEINO_p(1)\displaybreak[0]\\
  &~~ = \sum_{T_k\le t}\frac{2}{3} \frac{G^N(T_k)-G^N(T_{k-1})}{T_k-T_{k-1}}\left(\rho_{T_{k-1}}^2+1\right)\left(\sigma_{T_{k-1}}^X\sigma_{T_{k-1}}^Y\right)^2\Delta T_k +\KLEINO_p(1) \displaybreak[0]\\
  &~~\stackrel{p}{\longrightarrow} \frac{2}{3}\int_0^t(1+\rho_s^2)(\sigma_s^X\sigma_s^Y)^2G^{\prime}(s)\,ds~~.
\end{align*}
In the first step cross terms of the inner sums have been neglected since they are centred and by It\^{o} isometry it can be shown that their second moments are bounded from above by $C i^3\delta_N^3$. We frequently use estimates $\delta_N^{l-1}$ for sums of the type $\sum_i^N(\Delta T_i)^l\,,l>1,$ by Hölder's inequality with the supremum norm to obtain upper bounds.\\
Subsequently squared increments of $L$ and $M$ and the increments of the product $L\cdot M$ in these inner sums are substituted by the increments of the quadratic (co-)variation processes. The induced error terms are centred by It\^{o} isometry and involving Cauchy-Schwarz inequality it follows that $C\delta_N$
is an upper bound for their second moments. \\
The crucial non-standard approximation is that on each block $(T_{k-1},\ldots,T_{k-i\vee 0})$ the increments of the form $\int_{T_{k-l-1}}^{T_{k-l}}f(t)dt$ with continuous functions $f$ for $l=1,\ldots,k\wedge i$ are approximated by $\Delta T_k f\left(T_{k-1}\right)$. This blockwise approximation is treated in Lemma \ref{crucialh2} and makes use of the concept of a time-changed quadratic variation of times and particularly \eqref{tcaqvt}. Finally, $1/i\sum_{l=1}^{i}(1-(l/i))^2=1/3+\KLEINO(1)$ and the convergence in probability is ensured by Assumption \ref{grid3} and the convergence of the Riemann sums to the integral.
\begin{lem}\label{crucialh2}
On the same assumptions as in Proposition \ref{cltosdis}, it holds true that the term
\begin{align*}\frac{N}{iT}\sum_{T_k\le t}\left(\int_{T_{k-1}}^{T_k}(\sigma_s^X)^2ds\left(\sum_{l=1}^{i\wedge k}\left(1-\frac{l}{i}\right)^2\int_{T_{k-l-1}}^{T_{k-l}}(\sigma_s^Y)^2ds\right)-(\sigma_{T_{k-1}}^X\sigma_{T_{k-1}}^Y\Delta T_k)^2\sum_{l=1}^{i\wedge k}\left(1-\frac{l}{i}\right)^2\right)~,\end{align*}
and the analogous blockwise approximations of $\Delta [M]_{T_k}$ and $\Delta [M]_{T_k}$ through constant left-end points converge to zero in probability.
\end{lem}
\begin{proof}
The approximation uses the concept of a time-change in the asymptotic quadratic variation of refresh times introduced in \cite{zhang} which is expounded in \ref{tcr}. By virtue of that concept we may suppose without loss of generality that the sampling design of the closest synchronous approximation satisfies \eqref{tcaqvt}.\\
The asymptotic orders of the three terms are deduced analogously and we restrict us to the proof of the above given first term.
An application of the mean value theorem yields
\begin{align*}\frac{N}{iT}\sum_{T_k\le t}\int_{T_{k-1}}^{T_k}(\sigma_s^X)^2ds\left(\sum_{l=1}^{i\wedge k}\left(1-\frac{l}{i}\right)^2\int_{T_{k-l-1}}^{T_{k-l}}(\sigma_s^Y)^2ds\right)\\=\frac{N}{iT}\sum_{T_k\le t}(\sigma_{\zeta_k}^X)^2\Delta T_k\left(\sum_{l=1}^{i\wedge k}\left(1-\frac{l}{i}\right)^2(\sigma_{\zeta^*_{k-l}}^Y)^2\Delta T_{k-l}\right)\end{align*}
with $\zeta_k\in [T_{k-1},T_k]$, $\zeta^*_q\in [T_{q-1},T_q]$. Since the volatility processes $\sigma^X,\sigma^Y$ are uniformly continuous on $[0,T]$ by Assumption \ref{eff}
\begin{align*}\sum_{l=1}^{i\wedge k}\left(1-\frac{l}{i}\right)^2\left|(\sigma_{\zeta^*_{k-l}}^Y)^2-(\sigma_{T_{k-1}}^Y)^2\right|\hspace*{-.05cm}\Delta T_{k-l}\le i\delta_N\sup_{|t-s|\le i\delta_N}{\left|(\sigma_{t}^Y)^2-(\sigma_{s}^Y)^2\right|}=\KLEINO_{a.\,s.\,}\hspace*{-.1cm}(i\delta_N)~,\end{align*}
\begin{align*}\sum_{T_k\le t}\hspace*{-.05cm}\left|(\sigma_{\zeta_k}^X)^2\hspace*{-.05cm}-\hspace*{-.05cm}(\sigma_{T_{k-1}}^X)^2\right|(\sigma_{T_{k-1}}^Y)^2\Delta T_k\hspace*{-.05cm}\sum_{l=1}^{i\wedge k}\hspace*{-.05cm}\left(1-\frac{l}{i}\right)^2\hspace*{-.15cm}\Delta T_{k-l}\le i\delta_N\hspace*{-.15cm}\sup_{|t-s|\le \delta_N}\hspace*{-.075cm}{\left|(\sigma_{t}^Y)^2\hspace*{-.05cm}-\hspace*{-.05cm}(\sigma_{s}^Y)^2\right|}\hspace*{-.05cm}=\KLEINO_{a.\,s.\,}\hspace*{-.1cm}(i\delta_N\hspace*{-.05cm})\end{align*}
hold almost surely (denoted a.\,s.\,). \\
With the Cauchy-Schwarz inequality and \eqref{tcaqvt}, we obtain 
\begin{align*}\frac{N}{iT}\sum_{T_k\le t} &(\sigma_{T_{k-1}}^X\sigma_{T_{k-1}}^Y)^2\Delta T_k\left(\sum_{l=1}^{i\wedge k}\left(1-\frac{l}{i}\right)^2\left|\Delta T_{k-l}-\frac{T}{N}\right|\right)\\ &\le \frac{N}{iT}\sup_{s\in[0,t]}{\left(\sigma_t^X\sigma_t^Y\right)^2}\sum_{l=1}^{i}\sum_{j=1}^{N-l}\left|\left(\Delta T_j-\frac{T}{N}\right)\Delta T_{j+l}\right|\\ &\le \frac{N}{iT} C\left(\sum_{j=1}^N\left(T_{(j+i)\vee N}-T_j\right)^2\sum_{j=1}^N\left(\Delta T_j-\frac{T}{N}\right)^2\right)^{\nicefrac{1}{2}}=\KLEINO_{a.\,s.\,}(1)~.\end{align*}
Furthermore,
\begin{align*}\frac{N}{iT}\sum_{T_k\le t}(\sigma_{T_{k-1}}^X\sigma_{T_{k-1}}^Y)^2\Delta T_k\left|\Delta T_k-\frac{T}{N}\right|\sum_{l=1}^{i\wedge k}\left(1-\frac{l}{i}\right)^2&\hspace*{-.1cm}\le \frac{N}{T}C\hspace*{-.1cm}\left(\sum_{j=1}^N(\Delta T_j)^2\sum_{j=1}^N\left(\Delta T_j-\frac{T}{N}\right)^2\right)^{\nicefrac{1}{2}}\end{align*}
holds, where the right-hand side converges to zero almost surely due to \eqref{tcaqvt} and the Cauchy-Schwarz inequality. The preceding estimates imply the statement of the lemma.
\end{proof}
We proceed proving \eqref{covjac} that the quadratic covariations $\left[ \mathfrak{D}^N,L\right]_t$ and $\left[ \mathfrak{D}^N,M\right]_t$ converge to zero in probability for all $t\in[0,T]$.
$$\left[ \mathfrak{D}^{N},L\right]_t=\sqrt{\frac{N}{iT}}\sum_{T_k\le t}\hspace*{-.1cm}\left(\Delta\left[ L\right]_{T_k}\hspace*{-.1cm}\left(\sum_{l=1}^{i\wedge k}\left(1-\frac{l}{i}\right)\Delta M_{T_{k-l}}\right)\hspace*{-.1cm}+\Delta \left[ L,M\right]_{T_k}\left(\sum_{l=1}^{i\wedge k}\left(1-\frac{l}{i}\right)\Delta L_{T_{k-l}}\right)\right)$$
has an expectation equal to zero for all $t\in[0,T]$ and the second moment is bounded above by
$iN\delta_N^2$.
The order follows from the evaluation of the second moment that is carried out analogously as for the calculation of $\left[ \mathfrak{D}^{N}\right]_t$ before. For this reason $\left[ \mathfrak{D}^N,L\right]$ converges to zero in probability on $[0,T]$. It can be directly deduced that $\left[ \mathfrak{D}^N,M\right]_t=\KLEINO_p(1)$ as well. If $L^{\bot}$ is a bounded $(\mathcal{F}_t)$-martingale with $\left[ L,L^{\bot}\right]\equiv 0$, the quadratic covariation
$$\left[ \mathfrak{D}^N,L^{\bot}\right]_t=\sqrt\frac{N}{iT}\sum_{T_k\le t}\Delta\left[ L^{\bot},M\right]_{T_k}\left(\sum_{l=1}^{i\wedge k}\left(1-\frac{l}{i}\right)\Delta L_{T_{k-l}}\right)$$
converges to zero in probability on $[0,T]$ what can be concluded following the same principles and also that $\left[ \mathfrak{D}^N,M^{\bot}\right]_t=\KLEINO_p(1)~\forall M^{\bot}\in\mathcal{M}^{\bot},\,t\in[0,T]$. An application of Jacod's Theorem from \cite{jacod1} completes the proof of Proposition \ref{cltosdissynpr}.$\hfill\Box$\\ \noindent
\begin{prop}\label{cltmsdissynpr}On the same assumptions as in Proposition \ref{cltmsdis}, the continuous martingale 
\begin{align*}\mathfrak{M}_t^N\:=\sqrt{\frac{N}{M_N}}&\sum_{i=1}^{M_N}\left[\sum_{T_k\le t}\left(\Delta L_{T_k}+\int_{T_k}^t\sigma_s^XdW_s^X\right)\left(\sum_{l=1}^{i\wedge k}\left(1-\frac{l}{i}\right)\Delta M_{T_{k-l}}\right)\right. \\
&\left. +\sum_{T_k\le t}\left(\Delta M_{T_k}+\int_{T_k}^t\sigma_s^YdW_s^Y\right)\left(\sum_{l=1}^{i\wedge k} \left(1-\frac{l}{i}\right)\Delta L_{T_{k-l}}\right)\right]\end{align*}
for $t\in[0,T]$ converges stably in law as $N\rightarrow\infty, M_N\rightarrow\infty, M_N/N^{\alpha}\rightarrow 0$ for every $\alpha>2/3$ to a limiting time-changed Brownian motion
\begin{align*}\mathfrak{M}^N_t\stackrel{st}{\rightsquigarrow}\int_0^t\sqrt{v_{\mathfrak{M}_s}}d \tilde {\mathfrak{W}}_s^{\bot}~,\end{align*}
where $\tilde{\mathfrak{W}}^{\bot}$ is independent of $\mathcal{F}$ and with
\begin{align*} v_{\mathfrak{M}_s}=\frac{26}{35} T G^{\prime}(s)(\sigma_s^X\sigma_s^Y)^2(1+\rho_s^2)~.\end{align*}
\end{prop}
\noindent
\textit{Proof of Proposition \ref{cltmsdissynpr}:}\\
The discretization error of the generalized multiscale estimator calculated with the closest synchronous approximation under the equivalent martingale measure where the drift terms equal zero
\begin{align*}\sum_{i=1}^{M_N}\frac{\alpha_{i,M_N}^{opt}}{i}\sum_{j=i}^N\left(L_{T_j}-L_{T_{j-i}}\right)\left(M_{T_j}-M_{T_{j-i}}\right)-\left[ X,Y\right]_T\\
=\sum_{i=1}^{M_N}\alpha_{i,M_N}^{opt}\left(\frac{1}{i}\sum_{j=i}^N\left(L_{T_j}-L_{T_{j-i}}\right)\left(M_{T_j}-M_{T_{j-i}}\right)-\left[ X,Y\right]_T\right)\end{align*}
equals the weighted sum of $M_N\rightarrow\infty$ discretization errors of the type considered in Proposition \ref{cltosdis} because $\sum_{i=1}^{M_N}\alpha_{i,M_N}^{opt}=1$. Note, that all approximation errors in the preceding proof of Proposition \ref{cltosdissynpr} converge to zero in probability as long as $N\rightarrow\infty, i/N^{\alpha}\rightarrow 0$ for every $\alpha>2/3$.\\ 
We begin with the proof of a multivariate stable central limit theorem for a finite-dimensional vector:
\begin{lem}
\label{cltdisvec}
Consider the sequence of $K$-dimensional vectors $\mathds{D}^N=\left(\mathfrak{D}^{i_N^1},\ldots,\mathfrak{D}^{i_N^K}\right)$ where the entries ${\mathfrak{D}}^{i_N^k},k=1,\ldots,K<\infty$ are the continuous martingales
\begin{align*}{\mathfrak{D}}^{i_N^k}_t&=\left[\sum_{T_r\le t}\left(\Delta L_{T_r}+\int_{T_r}^t\sigma_s^XdW_s^X\right)\left(\sum_{l=1}^{i_N^k\wedge r}\left(1-\frac{l}{i_N^k}\right)\Delta M_{T_{r-l}}\right)\right. \\
&\left. +\sum_{T_r\le t}\left(\Delta M_{T_r}+\int_{T_r}^t\sigma_s^YdW_s^Y\right)\left(\sum_{l=1}^{i_N^k\wedge r} \left(1-\frac{l}{i_N^k}\right)\Delta L_{T_{r-l}}\right)\right]~\end{align*}
with a sequence of integers $i_N^k,k=1,\ldots,K$. On the Assumptions \ref{eff}, \ref{grid} and \ref{grid3} and if for every $k\in\{1,\ldots,K\}$ there exists a constant $q_k$ with $i_N^k/M_N\rightarrow q_k$, the following stable convergence holds true as $N\rightarrow\infty,M_N\rightarrow\infty,M_N/N^{\alpha}\rightarrow 0$ for every $\alpha>2/3$:
\begin{align}\sqrt{\frac{N}{M_N}}{\mathds{D}}^N_t\stackrel{st}{\rightsquigarrow}\int_0^t w_s d{\mathds{W}}_s~,\end{align}
with a $K$-dimensional Brownian motion $\mathds{W}$ independent of $\mathcal{F}$ and a predictable process $w_s$ with
\begin{align}\left(w_sw_s^*\right)_{mn}=\frac{T}{3}\min{(q_m,q_n)}\left(3-\frac{\min{(q_m,q_n)}}{\max{(q_m,q_n)}}\right)(1+\rho_s^2)\left(\sigma_s^X\sigma_s^Y\right)^2G^{\prime}(s)~\end{align}
with the convention that for $q_m=q_n=0$ the ratio is one.\\
For $\mathds{D}^N_T$ we obtain the following multivariate stable central limit theorem
\begin{align}\sqrt{\frac{N}{M_N}}{\mathds{D}}^N_T\stackrel{st}{\rightsquigarrow}\mathbf{N}\left(0,\eta^2\Sigma\right)~,\end{align}
with $\eta^2=2T\int_0^T(1+\rho_t^2)(\sigma_t^X\sigma_t^Y)^2G^{\prime}(t) dt$ and
$$\Sigma_{mn}=\frac{1}{6}\min{(q_m,q_n)}\left(3-\frac{\min{(q_m,q_n)}}{\max{(q_m,q_n)}}\right)~.$$
\end{lem}
\begin{proof}
Define for $k\in\{1,\ldots,K\}$ the continuous martingales
$$\mathfrak{M}_t^{i_N^k}=\sqrt{\frac{N}{M_N}}\mathfrak{D}_t^{i_N^k}~.$$
By virtue of Proposition \ref{cltosdis}, we already have that
$$\left[ \mathfrak{M}^{i_N^k}\right]_t\stackrel{p}{\longrightarrow}\frac{2}{3}T q_k\int_0^t(1+\rho_s^2)(\sigma_s^X\sigma_s^Y)^2G^{\prime}(s)\,ds~.$$
The limit of the quadratic covariations $\left[ \mathfrak{M}^{i_N^m},\mathfrak{M}^{i_N^k}\right]$ is derived using the same approximations as for the quadratic variation in the preceding proof:
{\allowdisplaybreaks[2]{
\begin{align*}\left[ \mathfrak{M}^{i_N^m},\mathfrak{M}^{i_N^k}\right]_t&=\frac{N}{M_N}\left[\sum_{T_r\le t}\Delta\left[ L\right]_{T_r}\left(\sum_{l=1}^{\min{(i_N^m,i_N^k,r)}}\left(1-\frac{l}{i_N^m}\right)\left(1-\frac{l}{i_N^k}\right)\left(\Delta M_{T_{r-l}}\right)^2\right)\right.\\
&\left. ~+ \sum_{T_r\le t}\Delta\left[ M\right]_{T_r}\left(\sum_{l=1}^{\min{(i_N^m,i_N^k,r)}}\left(1-\frac{l}{i_N^m}\right)\left(1-\frac{l}{i_N^k}\right)\left(\Delta L_{T_{r-l}}\right)^2\right)\right.\\
&\left. ~+\sum_{T_r\le t}2\Delta\left[ L,M\right]_{T_r}\left(\sum_{l=1}^{\min{(i_N^m,i_N^k,r)}}\left(1-\frac{l}{i_N^m}\right)\left(1-\frac{l}{i_N^k}\right)\Delta L_{T_{r-l}}\Delta M_{T_{r-l}}\right)\right]+\KLEINO_p(1)\\
&=N\sum_{T_r\le t}2\,\frac{G^{(N)}(T_r)-G^{(N)}(T_{r-1})}{\Delta T_r} (\rho_{T_{r-1}}^2+1)(\sigma_{T_{r-1}}^X\sigma_{T_{r-1}}^Y)^2\Delta T_r\\ &~~~~~~~~~~~~~~~~~~~~~~~~~~\times\left(\sum_{l=1}^{\min{(i_N^m,i_N^k,r)}}\left(1-\frac{l}{i_N^m}\right)\left(1-\frac{l}{i_N^k}\right)\right)+\KLEINO_p(1)\\
&\stackrel{p}{\longrightarrow} 2T\int_0^t(\rho_s^2+1)(\sigma_s^X\sigma_s^Y)^2G^{\prime}(s)\left(\frac{1}{6}\min{(q_m,q_k)}\left(3-\frac{\min{(q_m,q_k)}}{\max{(q_m,q_k)}}\right)\right)\,ds~,
\end{align*}}}
since $\sum_{l=1}^m(1-(l/m))(1-(l/M))=(1/2)m-(m^2/6M)-1/8+1/(12M)$ for $m,M\in\mathds{Z}$.
The multi-dimensional version of Jacod's stable central limit Theorem 2--1 from \cite{jacod1} enables us to prove the result of stable weak convergence of the vector provided we can verify the conditions
$$\left[ \mathds{D}^N,\mathds{L}\right]_t\stackrel{p}{\longrightarrow}0~~,~~\left[ \mathds{D}^N,\mathds{M}\right]_t\stackrel{p}{\longrightarrow}0,~~\forall t\in[0,T]~,$$
where $\mathds{L}$ denotes the vector with entries $\mathds{L}^{j}=L,j=1,\ldots,K$ and $\mathds{M}$ with $\mathds{M}^{j}=M,j=1,\ldots,K$, respectively,
and 
$$\left[ \mathds{D}^N,\mathds{L}^{\bot}\right]_t\stackrel{p}{\longrightarrow}0~~,~~\left[ \mathds{D}^N,\mathds{M}^{\bot}\right]_t\stackrel{p}{\longrightarrow}0,~~\forall t\in[0,T]$$
where $\mathds{L}^{\bot}$ and $\mathds{M}^{\bot}$ are bounded $(\mathcal{F}_t)$-adapted martingales orthogonal to $\mathds{L}$ and $\mathds{M}$, respectively.
That is because the reference continuous martingales for all entries of the vector $\mathds{D}^N$ are $L$ and $M$.
The componentwise proof of the conditions above is yet analogous as for the univariate case in the preceding proof. We conclude that the asymptotic distribution of the vector is described by a limiting time-changed Brownian motion on $[0,T]$, and the marginal distribution at time $T$ by a mixed Gaussian limit, where the normal distribution is defined as well as for all componentwise marginals on an orthogonal extension of the original underlying probability space.
\end{proof}
From the preceding multivariate limit theorem the Cram\'{e}r-Wold device allows to conclude the weak convergence of all one-dimensional linear combinations of the transformed discretization errors of a finite collection of one-scale subsampling estimators. For an asymptotically $\mathbf{N}\left(0,\Sigma\right)$-distributed random vector the sum of all components is asymptotically normally distributed with variance $\sum_{i,j}(\Sigma_{ij})$ by the Cram\'{e}r-Wold device and the normality of any linear sum of components of a multivariate normal distribution (see e.\,g.\,pp.\,516-517 in \cite{rao}).\\
The asymptotic variance in Proposition \ref{cltmsdis} is deduced from the multivariate limit and
$$\sum_{k,l}(\Sigma_{k,l})=2\sum_{k=1}^{M_N}\sum_{l=1}^k\frac{l}{6M_N}\left(3-\frac{l}{k}\right)\alpha_{k,M_N}^{opt}\alpha_{l,M_N}^{opt}+\KLEINO(1)=\frac{13}{35}+\KLEINO(1)$$
with the weights \eqref{optimalweights} inserted. \\
For the completion of the proof of Propositions \ref{cltmsdissynpr} and hence \ref{cltmsdis}, it remains to extend the result for asymptotically infinitely many addends. This part of the proof can be adopted from \cite{zhang} where a stable central limit theorem for a multiscale estimator for the integrated volatility in the univariate setting is proved. 
$\hfill\Box$
\subsubsection{Discretization error due to the lack of synchronicity}
\begin{prop}\label{Asy}On the Assumptions \ref{eff} and \ref{grid}, it holds true that
\begin{align*}\mathcal{A}_T^N&=\frac{1}{i}\sum_{j=i}^{N-1}\left[\left(L_{g_j}-L_{T_j}\right)\left(\sum_{k=j-i+1}^j\Delta M_{T_k}\right)+\left(M_{\gamma_j}-M_{T_{j}}\right)\left(\sum_{k=j-i+1}^j\Delta L_{T_k}\right)\right.\\
& \left.~~~~~~+ \Delta L_{T_{j+1}}\left(\sum_{k=j-i+1}^j\left(M_{T_k}-M_{\lambda_{k+1}}\right)\right)+\Delta M_{T_{j+1}}\left(\sum_{k=j-i+1}^j\left(L_{T_k}-L_{l_{k+1}}\right)\right)\right]=\mathcal{O}_p\left(\sqrt{i}{N}\right)~.\end{align*}
for the error associated with interpolation errors $\mathcal{A}_T^N$ for a one-scale subsampling estimator.
\end{prop}
\begin{proof}
$\mathcal{A}_T^N$ is the endpoint of a $\mathcal{F}_{j,N}=\mathcal{F}_{T_{j+1}^{(N)}}$-measurable discrete martingale with conditional expectation zero, since the addends incorporate products of Brownian increments over disjoint time intervals. The conditional variance yields
\begin{align*}
&\frac{1}{i^2}\sum_{j=i}^{N-1}\E\left[\left(\left(L_{g_j}-L_{T_j}\right)\left(\sum_{k=j-i+1}^j\Delta M_{T_k}\right)+\left(M_{\gamma_j}-M_{T_{j}}\right)\left(\sum_{k=j-i+1}^j\Delta L_{T_k}\right)\right.\right.\displaybreak[0]\\
& \left.\left.~~~~~~+ \Delta L_{T_{j+1}}\left(\sum_{k=j-i+1}^j\left(M_{T_k}-M_{\lambda_{k+1}}\right)\right)+\Delta M_{T_{j+1}}\left(\sum_{k=j-i+1}^j\left(L_{T_k}-L_{l_{k+1}}\right)\right)\right)^2\big|\mathcal{F}_{T_{j}^{(N)}}\right]\displaybreak[0]\\
&=\frac{1}{i^2}\sum_{j=i}^{N-1}\left(\E\left[(L_{g_j}-L_{T_j})^2\right]\left(\sum_{k=j-i+1}^j\Delta M_{T_k}\right)^2+\E\left[(M_{\gamma_j}-M_{T_{j}})^2\right]\left(\sum_{k=j-i+1}^j\Delta L_{T_k}\right)^2\right.\\
& \left.~~~~~~+ \E\left[(\Delta L_{T_{j+1}})^2\right]\left(\sum_{k=j-i+1}^j\left(M_{T_k}-M_{\lambda_{k+1}}\right)\right)^2+\E\left[(\Delta M_{T_{j+1}})^2\right]\left(\sum_{k=j-i+1}^j\left(L_{T_k}-L_{l_{k+1}}\right)\right)^2\right.\displaybreak[0]\\
&\left.~~~~~~+ \E\left[\int_{T_j}^{g_j}(\sigma_t^X)^2 dt\right]\left(\sum_{k=j-i+1}^j\left(M_{T_k}-M_{\lambda_{k+1}}\right)\right)\left(\sum_{k=j-i+1}^j\Delta M_{T_k}\right)\right.\displaybreak[0]\\
&\left. ~~~~~~+ \E\left[\int_{T_j}^{\gamma_j}(\sigma_t^Y)^2 dt\right]\left(\sum_{k=j-i+1}^j\left(L_{T_k}-L_{l_{k+1}}\right)\right)\left(\sum_{k=j-i+1}^j\Delta L_{T_k}\right)\right.\displaybreak[0]\\
&\left. ~~~~~~+ \E\left[\int_{T_j}^{g_j}\rho_t\sigma_t^X\sigma_t^Ydt\right]\left(\sum_{k=j-i+1}^j\left(L_{T_k}-L_{l_{k+1}}\right)\right)\left(\sum_{k=j-i+1}^j\Delta M_{T_k}\right)\right.\displaybreak[0]\\
&\left. ~~~~~~+ \E\left[\int_{T_j}^{\gamma_j}\rho_t\sigma_t^X\sigma_t^Ydt\right]\left(\sum_{k=j-i+1}^j\left(M_{T_k}-M_{\lambda_{k+1}}\right)\right)\left(\sum_{k=j-i+1}^j\Delta L_{T_k}\right)\right)\displaybreak[0]\\
&=\mathcal{O}_p\left(i^{-1}N^{-1}\right)~.
\end{align*}
The variance of the term is of order $(iN)^{-1}$ which can be proved by taking the expectation of the above given conditional variance and an upper bound of the second moment. The asymptotic orders of the addends follow from taking the expectations using It\^{o} isometry and analyzing the differences of the addends minus their expectations, that converge to zero at a faster rate. That part is similar to the proofs above and we forgo a more detailed computation here.
\end{proof}
Denote $A_T^{N,i}$ the error due to non-synchronicity and interpolations for a fixed subsampling frequency $i=1,\ldots,M_N$ in the following. The error due to asynchronicity of the generalized multiscale estimator \eqref{multiscale} equals the weighted sum
$\sum_{i=1}^{M_N}\alpha_{i,M_N}^{opt}A_T^{N,i}$. It has expectation zero and the variance is of order
\begin{align*}\var\left(\sum_{i=1}^{M_N}\alpha_{i,M_N}^{opt}\mathcal{A}_T^{N,i}\right)&=\sum_{i,k}\alpha_{i,M_N}^{opt}\alpha_{k,M_N}^{opt}\cov\left(\mathcal{A}_T^{N,i}\,,\,\mathcal{A}_T^{N,k}\right)\\ &=\underbrace{\sum_{i=1}^{M_N}\left(\alpha_{i,M_N}^{opt}\right)^2\E\left[\left(\mathcal{A}_T^{N,i}\right)^2\right]}_{=\mathcal{O}\left(M_N^{-2}N^{-1}\right)}+\underbrace{\sum_{i\ne k}\alpha_{i,M_N}^{opt}\alpha_{k,M_N}^{opt}\E\left[\mathcal{A}_T^{N,i}\mathcal{A}_T^{k,N}\right]}_{=\mathcal{O}\left(M_N^{-1}N^{-1}\right)}=\KLEINO\left(\frac{M_N}{N}\right)~.\end{align*}
Thus, the error due to interpolations is of smaller asymptotic order than the discretization error of the closest synchronous approximation and asymptotically negligible.
\subsection{Asymptotics of the cross term}
For a one-scale subsampling estimator cross terms are asymptotically negligible and hence the stable central limit theorem in Theorem \ref{cltone} is implied by Theorem \ref{cltosdis}. For the proof of the stable central limit theorem in Theorem \ref{cltmulti} for the multiscale approach, we cope with the asymptotics of the cross terms in this subsection. 
\begin{prop}\label{cltcross}On the Assumptions \ref{eff}, \ref{grid}, \ref{grid4} and \ref{e}, the cross terms of the generalized multiscale estimator \eqref{multiscale} with noise-optimal weights \eqref{optimalweights} weakly converge to a mixed normal limit as $M_N\rightarrow\infty,\,N\rightarrow\infty,\,M_N \delta_N\rightarrow 0$:
\begin{align}\notag\sqrt{M_N}\sum_{i=1}^{M_N}\frac{\alpha_{i,M_N}^{opt}}{i}\sum_{j=i}^N\left((X_{g_j}-X_{l_{j-i+1}})(\epsilon^Y_{\gamma_j}-\epsilon^Y_{\lambda_{j-i+1}})+(Y_{\gamma_j}-Y_{\lambda_{j-i+1}})(\epsilon^X_{g_j}-\epsilon^X_{l_{j-i+1}})\right)\\ \rightsquigarrow\mathbf{N}\left(0,\AVAR_{cross}\right)~,\end{align}
with asymptotic variance
\begin{align}\label{avarcross}\AVAR_{cross}=\frac{12}{5}\left(\eta_Y^2\int_0^T(1+I^{\prime}_Y(t))(\sigma_t^X)^2\,dt\,+\eta_X^2\int_0^T(1+I^{\prime}_X(t))(\sigma_t^Y)^2\,dt\right)~.\end{align}
The convergence holds conditionally given the paths of the efficient processes.
\end{prop}
\begin{proof}
This proof  affiliates to the discussion in Section \ref{sec:4}, where degrees of regularity of non-synchronous sampling schemes have been defined in Definition \ref{dr} that are assumed to converge to continuously differentiable functions.\\
On the Assumption \ref{e} of independent observation noise of $X$ and $Y$, the two different cross terms are uncorrelated and we prove a central limit theorem for the first one:
$$\sqrt{M_N}\sum_{i=1}^{M_N}\frac{\alpha_{i,M_N}^{opt}}{i}\sum_{j=i}^N(X_{g_j}-X_{l_{j-i+1}})(\epsilon^Y_{\gamma_j}-\epsilon^Y_{\lambda_{j-i+1}})\rightsquigarrow\mathbf{N}\left(0,\frac{12}{5}\eta_Y^2\int_0^T(1+I^{\prime}_Y(t))(\sigma_t^X)^2\,dt\right)~.$$
The parallel result for the other term can be proved analogously.\\
For the purpose of a shorter notation we have left out superscripts of the observation times, and write $\alpha_i,\,i=1,\ldots,M_N$ for the weights although we are interested in the specific weights \eqref{optimalweights}. Denote $\delta_N=\sup_{i\in\{1,\ldots,N\}}\Delta T_i$ and $\gamma_{j,+}=\min{\left(\tau_k\in \mathcal{T}^Y|\tau_k\in\G^{j+1}\right)},g_{j,+}=\min{\left(t_k\in \mathcal{T}^X|t_k\in\H^{j+1}\right)}$ and $C$ a generic constant as before.
From
\begin{align*}&\E\left[\left(\sqrt{M_N}\sum_{i=1}^{M_N}\frac{\alpha_i}{i}\sum_{j=i}^N(X_{g_j}-X_{T_j})(\epsilon^Y_{\gamma_j}-\epsilon^Y_{\lambda_{j-i+1}})+(X_{T_{j-i}}-X_{l_{j-i+1}})(\epsilon^Y_{\gamma_j}-\epsilon^Y_{\lambda_{j-i+1}})\right)^2\right]\\
&~~\le M_N\sum_{i,k\in\{1,\ldots,M_N\}}\frac{\alpha_i\alpha_k}{ik}2\eta_Y^2\left(\sum_{j=i\vee k}^N\E(X_{g_j}-X_{T_j})^2+\sum_{j=0}^{N-(i\vee k)}\E(X_{T_j}-X_{l_{j+1}})^2\right)\\
&~~\le M_N\, C 4\eta_Y^2\sum_{i,k\in\{1,\ldots,M_N\}}\frac{\alpha_i\alpha_k}{ik}=\mathcal{O}\left(M_N^{-1}\right)~,
\end{align*}
for the errors due to interpolations and 
\begin{align*}&\E\left[\left(\sqrt{M_N}\sum_{i=1}^{M_N}\frac{\alpha_i}{i}\left(\sum_{k=N-i+1}^N\epsilon^Y_{\gamma_k}(X_{T_k}-X_{T_{k-i}})-\sum_{k=1}^{i}\epsilon^Y_{\lambda_k}(X_{T_{k+i}}-X_{T_k})\right)\right)^2\right]\\
&~~=M_N\sum_{i,k\in\{1,\ldots,M_N\}}\frac{\alpha_i\alpha_k}{ik}\eta_Y^2\left(\sum_{r=N-(i\wedge k)+1}^N\E(X_{T_r}-X_{T_{r-i}})^2+\sum_{r=1}^{i\wedge k}\E(X_{T_{r+i}}-X_{T_r})^2\right)\\ &~~=\mathcal{O}\left(M_N \delta_N\right)
\end{align*}
for boundary terms, we conclude that
\begin{align*}&\sqrt{M_N}\sum_{i=1}^{M_N}\frac{\alpha_i}{i}\sum_{j=i}^N(X_{g_j}-X_{l_{j-i+1}})(\epsilon^Y_{\gamma_j}-\epsilon^Y_{\lambda_{j-i+1}})\\
&=\sqrt{M_N}\sum_{i=1}^{M_N}\frac{\alpha_i}{i}\left(\sum_{j=i}^{N-i}\epsilon^Y_{\gamma_j}(X_{T_j}-X_{T_{j-i}})-\epsilon^Y_{\lambda_{j+1}}(X_{T_{j+i}}-X_{T_j})\right)+\KLEINO_p(1)\\
&=\sqrt{M_N}\sum_{j=2}^{N-2}\left(\epsilon^Y_{\gamma_j}\sum_{i=1}^{M_N^*(j)}\frac{\alpha_i}{i}(X_{T_j}-X_{T_{j-i}})-\epsilon^Y_{\lambda_{j+1}}\sum_{i=1}^{M_N^*(j)}(X_{T_{j+i}}-X_{T_j})\right)+\KLEINO_p(1)\\	
&=\sqrt{M_N}\left(\sum_{j\in\mathcal{Y}_1}\epsilon^Y_{\gamma_j}\sum_{i=1}^{M_N^*(j)}\frac{\alpha_i}{i}\zeta_{i,j}^1+\sum_{j\in\mathcal{Y}_2}\epsilon^Y_{\gamma_j}\sum_{i=1}^{M_N^*(j)}\frac{\alpha_i}{i}\zeta_{i,j}^2 +\sum_{j\in\mathcal{Y}_3}\epsilon^Y_{\gamma_j}\sum_{i=1}^{M_N^*(j)}\frac{\alpha_i}{i}\zeta_{i,j}^3\right. \\
& \left.~~~~~~~~~~~~~~~~~~~~~+\sum_{j\in\mathcal{Y}_4}\epsilon^Y_{\gamma_j}\sum_{i=1}^{M_N^*(j)}\frac{\alpha_i}{i}\zeta_{i,j}^{4a}-\sum_{j\in\mathcal{Y}_4}\epsilon^Y_{\gamma_{j,+}}\sum_{i=1}^{M_N^*(j)}\frac{\alpha_i}{i}\zeta_{i,j}^{4b}\right)+\KLEINO_p(1)~.
\end{align*}
Here, we aggregate the observation times $\gamma_j,\lambda_j,~j=2,\ldots,N-2$ in disjoint sets conforming to the four cases discussed in Section \ref{sec:4}. Denote thereto
\begin{align*}\mathcal{Y}_1&=\{j\in\{2,\ldots,N_2\}|\gamma_j \ne \gamma_{j-1}\,,\,\gamma_j\le g_j\}~,\\
\mathcal{Y}_2&=\{j\in\{2,\ldots,N_2\}|\gamma_j > g_j\,,\,\gamma_j\ge g_{j,+}\}~,\\
\mathcal{Y}_3&=\{j\in\{2,\ldots,N_2\}|\gamma_j > g_j\,,\,\gamma_j< g_{j,+}\,,\,\gamma_{j,+}>g_{j,+}\}~,\\
\mathcal{Y}_4&=\{j\in\{2,\ldots,N_2\}|\gamma_j > g_j\,,\,\gamma_j< g_{j,+}\,,\,\gamma_{j,+}\le g_{j,+}\}~,\end{align*}
and $M_N^*(j)=\min{(j,N-j,M_N)}$.
The increments of $X$ that are multiplied with each observation error differ according to the set $\mathcal{Y}_k,~1\le k\le 4$ to which $\gamma_j$ belongs. We use the notation
\begin{align*}\zeta_{i,j}^1&=(X_{T_j}-X_{T_{j-i}})-(X_{T_{j+i}}-X_{T_j})~,\\
							\zeta_{i,j}^2&=(X_{T_j}-X_{T_{j-i}})+(X_{T_{j+1}}-X_{T_{j-i+1}})-(X_{T_{j+i+1}}-X_{T_{j+1}})~,\\
							\zeta_{i,j}^3&=(X_{T_j}-X_{T_{j-i}})-(X_{T_{j+i+1}}-X_{T_{j+1}})~,\\
							\zeta_{i,j}^{4a}&=(X_{T_j}-X_{T_{j-i}})~,~\zeta_{i,j}^{4b}=(X_{T_{j+i+1}}-X_{T_{j+1}})~.\end{align*}
							
The resulting aggregated leading term above of the cross term is the endpoint of a discrete martingale with respect to the filtration $\mathcal{F}_{j,N}\:=\sigma\big(\epsilon^Y_{\tau_k}|	\tau_k<\gamma_{j+1}\,,\,X,Y	\big)$. Since if $j\in\mathcal{Y}_4\Rightarrow \gamma_{j,+}<\gamma_{j+1}$, the martingale property with respect to the filtration $\mathcal{F}_{j,N}$	is assured by Assumption \ref{e}.		\\
An application of the non-stable version of the central limit theorem for martingale triangular arrays from \cite{hall} will proof the asymptotic normality of the cross term conditionally on the paths of the efficient processes. The conditional Lindeberg condition can be verified (using Chebyshev's inequality or directly verifying the conditional Lyapunov condition) in the same way as before and we omit it here. The sum of conditional variances yields
\begin{align*}&\sum_{l\in\{1,2,3,4a\}}\hspace*{-0.15cm}\left(\sum_{j\in\mathcal{Y}_l}\hspace*{-0.075cm}\E\hspace*{-0.075cm}\Big[\Big(\sqrt{M_N}\epsilon^Y_{\gamma_j}\sum_{i=1}^{M_N^*(j)}\frac{\alpha_i}{i}\zeta_{i,j}^l\Big)^2\Big| \mathcal{F}_{j-1,N}\hspace*{-0.075cm}\Big]\hspace*{-0.075cm}+\hspace*{-0.1cm}\sum_{j\in\mathcal{Y}_4}\hspace*{-0.075cm}\E\hspace*{-0.075cm}\Big[\Big(\hspace*{-0.125cm}-\hspace*{-0.075cm}\sqrt{M_N}\epsilon^Y_{\gamma_{j,+}}\sum_{i=1}^{M_N^*(j)}\frac{\alpha_i}{i}\zeta_{i,j}^{4b}\Big)^2\Big| \mathcal{F}_{j-1,N}\hspace*{-0.05cm}\Big]\hspace*{-0.05cm}\right)\\
&=M_N\eta_Y^2\left(\sum_{j\in\mathcal{Y}_1\cup\mathcal{Y}_3\cup\mathcal{Y}_4}\left(\sum_{i=1}^{M_N^*(j)}\frac{\alpha_i}{i}\zeta_{i,j}^1\right)^2+\sum_{j\in\mathcal{Y}_2}\left(\sum_{i=1}^{M_N^*(j)}\frac{\alpha_i}{i}\zeta_{i,j}^2\right)^2\right)+\KLEINO_p(1)\\
&=M_N\eta_Y^2\left(\sum_{j\in\mathcal{Y}_1\cup\mathcal{Y}_3\cup\mathcal{Y}_4}\sum_{i,k\in\{1,\ldots,M_N^*(j)\}}\hspace{-.25cm}\frac{\alpha_i\alpha_k}{ik}\left(\zeta_{i\wedge k,j}^1\right)^2\hspace{-.1cm}+\hspace{-.1cm}\sum_{j\in\mathcal{Y}_2}\hspace{-.1cm}\left(\sum_{i,k\in\{1,\ldots,M_N^*(j)\}}\frac{\alpha_i\alpha_k}{ik}\left(\zeta_{i\wedge k,j}^2\right)^2\right)\hspace{-.1cm}\right)\hspace{-.1cm}+\KLEINO_p(1)\\
&=M_N\eta_Y^2\left(\sum_{j\in\mathcal{Y}_1\cup\mathcal{Y}_3\cup\mathcal{Y}_4}\sum_{i,k\in\{1,\ldots,M_N^*(j)\}}\frac{\alpha_i\alpha_k}{ik}\left((X_{T_j}-X_{T_{j-(i\wedge k)}})^2+(X_{T_{j+(i\wedge k)}}-X_{T_j})^2\right)\right.\\
&\left. ~~~~~~~~~~~~~~+\sum_{j\in\mathcal{Y}_2}\sum_{i,k\in\{1,\ldots,M_N^*(j)\}}\frac{\alpha_i\alpha_k}{ik}\left(4(X_{T_j}-X_{T_{j-(i\wedge k)}})^2+(X_{T_{j+(i\wedge k)}}-X_{T_j})^2\right)\right)+\KLEINO_p(1)\\
&=M_N\eta_Y^2\sum_{j=2}^{N-2}\sum_{i,k\in\{1,\ldots,M_N^*(j)\}}\frac{\alpha_i\alpha_k}{ik}(2+\1_{\{j\in\mathcal{Y}_2\}})(X_{T_j}-X_{T_{j-(i\wedge k)}})^2+\KLEINO_p(1)\\
&=M_N\eta_Y^2\left(\sum_{i,k\in\{1,\ldots,M_N^*(j)\}}\frac{\alpha_i\alpha_k}{ik}\left(2(i\wedge k)\widehat{\left[ X\right]}_T^{sub,i\wedge k}+\sum_{j=i\wedge k}^N\1_{\{j\in\mathcal{Y}_2\}}(X_{T_j}-X_{T_{j-(i\wedge k)}})^2\right)\right)+\KLEINO_p(1)\\
&\stackrel{p}{\longrightarrow}\frac{12}{5}\eta_Y^2\left(\left[ X\right]_T+\int_0^TI^{\prime}_Y(t)(\sigma_t^X)^2\,dt\right)~.
\end{align*}
Since for the shifted increments 
$$(X_{T_{j+i+1}}-X_{T_{j+1}})=(X_{T_{j+i}}-X_{T_{j}})+\mathcal{O}_p\left(N^{-\nicefrac{1}{2}}\right)$$
holds, where the order is for time instants of average length $N^{-1}$, the variances of the sums over all $j\in\mathcal{Y}_1$ and $j\in\mathcal{Y}_3$ are asymptotically equal. The variance of both uncorrelated sums over maxima $\gamma_j$ and minima $\gamma_{j,+}$ distributed according to the fourth case is also asymptotically equal to the variances of those two addends. Only the asymptotic variance of the sum over all $j\in\mathcal{Y}_2$ is bigger. For this reason the total asymptotic variance hinges on the asymptotic degree of regularity of the non-synchronous sampling scheme $(\mathcal{T}^X,\mathcal{T}^Y)$ defined in Definition \ref{dr}.\\
In the calculation of the asymptotic variance we have used that
$$\zeta_{i,j}^1\zeta_{i,k}^1=\left(\zeta_{i\wedge k,j}^1\right)^2+\zeta^1_{i\wedge k,j}\left(\sum_{l=j-(i\vee k)+1}^{j-(i\wedge k)}\Delta X_{T_l}+\sum_{l=j+(i\wedge k)+1}^{j+(i\vee k)}\Delta X_{T_l}\right)~,$$
where the second remainder addend has an expectation equal to zero, and analogous formulae for $\zeta_{i,j}^2$, for all $1\le i\le M_N\,,1\le k\le M_N\,,\,k\vee i\le j\le N-(i\vee k)$.\\
Furthermore, an application of the mean value theorem, It\^{o} isometry and approximations in the same spirit as in the calculation of the asymptotic variance in the proof of the central limit theorem for the discretization errors of the estimators, lead to the Riemann sum in the calculation of the asymptotic variance above. The cross terms in $(\zeta_{i,j}^l)^2$, $l=1,2$ are asymptotically negligible. Since in $\mathcal{Y}_4$ repeating maxima $\gamma_i=\gamma_{i+1}$ are considered only once, it holds true that $|\mathcal{Y}_1|+|\mathcal{Y}_3|+|\mathcal{Y}_4|+2|\mathcal{Y}_2|=N-3\pm 1$ (the last addend can appear due to boundary term effects). In the last step we have used that
$$M_N\sum_{i,k\in\{1,\ldots,M_N\}}\frac{\alpha_{i,M_N}^{opt}\alpha_{k,M_N}^{opt}}{ik}(i\wedge k)=6/5+\KLEINO(1)$$
when inserting the weights \eqref{optimalweights}.\\
From the analysis for the asymptotic discretization error of a one-scale subsampling estimator, we know that
$$\widehat{\left[ X\right]}_T^{sub,i\wedge k}=\frac{1}{i\wedge k}\sum_{j=i\wedge k}^N(X_{T_j}-X_{T_{j-(i\wedge k)}})^2=\left[ X\right]_T+\mathcal{O}_p\left(\sqrt{\frac{(i\wedge k)}{N}}\right)$$
holds true. Similarly, it can be deduced that
\begin{align*}\frac{1}{i}\sum_{j=i}^N\1_{\{j\in\mathcal{Y}_2\}}(X_{T_j}-X_{T_{j-i}})^2&=\frac{1}{i}\sum_{l=1}^N\left(\Delta X_{T_l}\right)^2\left(\sum_{k=1}^{i}\1_{\{(k+l-1)\in\mathcal{Y}_2\}}\right)+\mathcal{O}_p\left(\sqrt{i/N}\right)\\
&=\int_0^T I^{\prime}_Y(t)(\sigma_t^X)^2\,dt+\mathcal{O}_p\left(\sqrt{i/N}\right)~,
\end{align*}
on Assumption \ref{grid4}. 
\end{proof}

\section{Proof of Proposition \ref{avarest}\label{sec:8}}
Let $R_j^N$, $R_j^n$, $R_j^m$, $S_j^{N,X}$ and $S_j^{N,Y}$ denote the number of times $T_k^{(N)},0\le k\le N$, $t_i^{(n)}$ and $\tau_j^{(m)}$ in the bins $[G_j^N,G_{j+1}^N)$, $[(I_Y)_j^N,(I_Y)_{j+1}^N)$, $[(I_X)_j^N,(I_X)_{j+1}^N) ,0\le j\le K_N-1,$. Define the generalized multiscale estimator in the fashion of \eqref{multiscale}
$$\widehat{\Delta\left[ X,Y\right]}_{G_{j+1}^N}=\sum_{i=1}^{M_N(j)}\frac{\alpha_{i,M_N(j)}^{opt}}{i}\sum_{r=i}^{R_j^N}\left(\tilde X_{g_r}-\tilde X_{l_{r-i+1}}\right)\left(\tilde Y_{\gamma_r}-\tilde Y_{\lambda_{r-i+1}}\right)$$
for the increase of the quadratic covariation $\Delta\left[ X,Y\right]_{G_{j+1}^N}$ and the univariate multiscale estimators
\begin{align*}\widehat{\Delta\left[ X\right]}_{G_{j+1}^N}=\sum_{i=1}^{M_n(j)}\frac{\alpha_{i,M_n(j)}^{opt}}{i}\sum_{r=i}^{R_j^n}\left(\tilde X_{t_r}-\tilde X_{t_{r-i}}\right)^2,\;
\widehat{\Delta\left[ Y\right]}_{G_{j+1}^N}=\sum_{i=1}^{M_m(j)}\frac{\alpha_{i,M_m(j)}^{opt}}{i}\sum_{r=i}^{R_j^m}\left(\tilde Y_{\tau_r}-\tilde Y_{\tau_{r-i}}\right)^2~,\end{align*}
and $\widehat{\Delta\left[ X\right]}_{(I_Y)_{j+1}^N}, \widehat{\Delta\left[ Y\right]}_{(I_X)_{j+1}^N}$ analogously, where all binwise multiscale frequencies are of order $\sqrt{N K_N}$ with possibly differing constants. Essential when considering the multiscale estimators on bins is that on Assumption \ref{grid3} the distances between sampling times are of order $N^{-1}\sim n^{-1}\sim m^{-1}$, whereas the numbers of observations $R_j^{\cdot},S_j^{\cdot},1\le j\le K_N$ in the specific bin are at most of order $NK_N^{-1}$. Following the analysis for the four uncorrelated parts of the estimation error in sections \ref{sec:7}, orders of the discretization variances 
$$\sum_{i,k\in\{1,\ldots,M_N(j)\}}\frac{\alpha_{i,M_N(j)}^{opt}\alpha_{k,M_N(j)}^{opt}}{ik}\cdot i\cdot R_j^N\frac{i^2}{N^2}\sim M_N(j)\frac{R_j^N}{N^2}\sim\frac{M_N(j)}{K_NN}~,$$
$M_n(j)/(nK_N)$, $M_m(j)/(mK_N)$,$M_n(j)\frac{S_j^{N,X}}{N^2}$ and $M_n(j)\frac{S_j^{N,Y}}{N^2}$, respectively.
Cross terms are of order $R_j^N/(NM_N(j))\sim (M_N(j)K_N)^{-1}$ and analogous orders for the univariate estimators.
The errors due to noise instead depend only on the number of observations in the considered interval. Therefore, the addends are of orders $R_j^N/M_N^3(j)\sim\frac{N}{K_NM_N^3(j)}$ and $M_N^{-1}(j)$ and analogous for the univariate estimators.\\
Choosing all multiscale frequencies $M_{\cdot}(j)\sim N^{\nicefrac{1}{2}}K_N^{\nicefrac{1}{2}}$ for every $j$, so that $M_N(\cdot)N^{\nicefrac{1}{2}}\rightarrow\infty$, the error due to end-effects in the noise part and the discretization error dominate asymptotically the two other addends and are of order $N^{-\nicefrac{1}{4}}K_N^{-\nicefrac{1}{4}}$. This holds as long as $K_N N^{-\nicefrac{1}{3}}\rightarrow 0$, such that $M_N(j)(N/K_N)^{-1}\rightarrow 0$ as $N\rightarrow\infty$.\\
The estimators \eqref{avarest1}-\eqref{avarest4} are consistent as $K_N\rightarrow\infty$ with $K_NN^{-\nicefrac{1}{3}}\rightarrow 0$ as $N\rightarrow \infty$, since
\begin{align*}\hat I_2&=\sum_{j=1}^{K_N}\left(\frac{\widehat{\Delta\left[ X\right]}_{G_j^N}\widehat{\Delta\left[ Y\right]}_{G_j^N}}{\left(\Delta G_j^N\right)^2}\right)\hspace*{-.05cm}\frac{G^N(T)}{K_N}\\
											&=\sum_{j=1}^{K_N}\left(\frac{\int_{G_{j-1}^N}^{G_j^N}\left(\sigma_t^X\right)^2\,dt\int_{G_{j-1}^N}^{G_j^N}\left(\sigma_t^Y\right)^2\,dt+\mathcal{O}_p\left(N^{-\nicefrac{1}{4}}K_N^{-\nicefrac{1}{4}}\right)}{\Delta G_j^N}\right)^2\,\frac{G^N(T)}{K_N}\\
											&=\sum_{j=1}^{K_N}\hspace*{-.05cm}\left(\sigma^X\right)_{\overline{G_j^N}}^2\left(\sigma^Y\right)_{\widetilde{G_j^N}}^2\frac{G^N(T)}{K_N}+\hspace*{-.05cm}\mathcal{O}_p\left(K_N^{\nicefrac{1}{4}}N^{-\nicefrac{1}{4}}\right)
											=\sum_{j=1}^{K_N}\hspace*{-.05cm}\left(\sigma^X\sigma^Y\right)_{G_{j-1}^N}^2\frac{G^N(T)}{K_N}+\hspace*{-.05cm}\mathcal{O}_p\left(K_N^{\nicefrac{1}{4}}N^{-\nicefrac{1}{4}}\right)
											\end{align*}
and similar conclusions for the other three estimators. We have used that $\Delta G_j^N\sim N^{-1}$	and apply the mean value theorem. $\overline{G_j^N}$ is some value $G_{j-1}^N\le \overline{G_j^N}\le G_j^N$. Finally, elementary inequalities as
{\small
$\left|\left(\sigma^X\right)^2_{\overline{G_j^N}}\left(\sigma^Y\right)^2_{\widetilde{G_j^N}}-\left(\sigma^X\sigma^Y\right)^2_{G_{j-1}^N}\right|~\le \left(\sigma^X\right)^2_{\overline{G_j^N}}\left|\left(\sigma^Y\right)^2_{\widetilde{G_j^N}}-\left(\sigma^Y\right)^2_{G_{j-1}^N}\right|+\left(\sigma^Y\right)^2_{G_{j-1}^N}\left|\left(\sigma^X\right)^2_{\overline{G_j^N}}-\left(\sigma^X\right)^2_{G_{j-1}^N}\right|$}\normalsize
and 
{\small $\sum_{j=1}^{K_N}\left|\left(\rho\sigma^X\sigma^Y\right)_{\overline{G_j^N}}^2-\left(\rho\sigma^X\sigma^Y\right)_{G_{j-1}^N}^2\right|\,\frac{G^N(T)}{K_N}\le \sup_{|t-s|\le \Delta \sup_j{G_j^n}}{\left|\rho_t\sigma_t^X\sigma_t^Y-\rho_s\sigma_s^X\sigma_s^Y\right|}G^N(T)=\KLEINO_{a.\,s.\,}(1)$}\normalsize
are involved in the approximations of the type above.
Considering \eqref{avarest3} and \eqref{avarest4}, note that bin-widths chosen accordingly to $I_Y^N$ are asymptotically of order $K_N^{-1}$ in any interval of $[0,T]$ on that the corresponding part of the integral $\int I_Y^{\prime}(t)(\sigma_t^X)^2\,dt$ is strictly positive. \\
Denote $R_N^k,\,k=1,\ldots,4$, the orders of the approximation errors of the four above given integrals and their Riemann sums evaluated on the partition given $K_N$ bins. The variance of the estimators $\widehat{\eta_X^2}$ and $\widehat{\eta_Y^2}$ for the noise variances are known to be $\E\left[\left(\epsilon_{t_1}^X\right)^4\right]N^{-1}$ and $\E\left[\left(\epsilon_{\tau_1}^Y\right)^4\right]N^{-1}$ and hence $\mathcal{O}\left(N^{-1}\right)$ on Assumption \ref{e} from \cite{zhangmykland}. From 
\begin{align*}\hat I_k=I_k+\mathcal{O}_p\left(R_N^k+K_N^{\nicefrac{1}{2}}N^{-\nicefrac{1}{2}}\right)~,~k=1,\ldots,4~,\end{align*}
we derive that
$$\widehat{\AVAR}_{multi}=\AVAR_{multi}+\mathcal{O}_p\left(\max_k{R_N^k}+K_N^{\nicefrac{1}{2}} N^{-\nicefrac{1}{2}}\right)~,$$
and the same result for the one-scale estimator. $\hfill\Box$\\ \noindent
\bibliographystyle{model1b-num-names}
\bibliography{literatur}					

\end{document}